\documentclass[a4paper,10pt]{article}
\usepackage{amsmath,amsthm,latexsym,makeidx,amscd,amssymb}
\usepackage[all]{xy}
\usepackage[hypertex,backref]{hyperref}
\usepackage{graphics}
\makeindex
\numberwithin{equation}{subsection}

\newtheorem{prop}{Proposition}[section]
\newtheorem{lem}[prop]{Lemma}
\newtheorem{ddd}[prop]{Definition}
\newtheorem{theorem}[prop]{Theorem}
\newtheorem{cor}[prop]{Corollary}

\newcommand{\F}{{\mathcal F}}
\newcommand{\Pj}{{\mathcal P}}

\newcommand{\Aut}{\mathop{\mbox{\rm Aut}}}

\newcommand{\Ch}{\mbox{\rm Ch}}
\newcommand{\End}{\mathop{\mbox{\rm End}}}
\newcommand{\Td}{\mathop{\mbox{\rm Td}}}

\newcommand{\ind}{\mathop{\mbox{\rm ind}}}
\newcommand{\Ind}{{\rm ind}}

\newcommand{\dom}{\mathop{\rm dom}}

\newcommand{\Coker}{\mathop{\rm Coker}}

\newcommand{\tr}{\mathop{\rm tr}}

\newcommand{\diag}{\mathop{\rm diag}}
\newcommand{\chh}{{\rm ch}}

\newcommand{\A}{{\mathcal A}}
\newcommand{\B}{{\mathcal B}}
\newcommand{\Ai}{{\mathcal A}_{\infty}}
\newcommand{\Bi}{{\mathcal B}_{\infty}}

\newcommand{\C}{C^{\infty}}

\newcommand{\Ok}{\hat \Omega_k}
\newcommand{\Oi}{\hat \Omega_*}
\newcommand{\Ob}{\hat \Omega}
\newcommand{\ra}{\partial}
\newcommand{\lr}{\longrightarrow}

\newcommand{\ten}{\otimes}

\newcommand{\ov}{\overline}

\newcommand{\incl}{\hookrightarrow}

\newcommand{\dira}{\partial \!\!\!/}

\DeclareMathOperator{\ch}{ch}

\DeclareMathOperator{\supp}{supp}
\DeclareMathOperator{\spfl}{sf}

\DeclareMathOperator{\Ker}{Ker}
\DeclareMathOperator{\di}{d}

\def\bbbr{{\rm I\!R}} 
\def\bbbn{{\rm I\!N}} 

\def\bbbc{{\rm I\!C}}

\def\bbbq{{\mathchoice {\setbox0=\hbox{$\displaystyle\rm Q$}\hbox{\raise
0.15\ht0\hbox to0pt{\kern0.4\wd0\vrule height0.8\ht0\hss}\box0}}
{\setbox0=\hbox{$\textstyle\rm Q$}\hbox{\raise
0.15\ht0\hbox to0pt{\kern0.4\wd0\vrule height0.8\ht0\hss}\box0}}
{\setbox0=\hbox{$\scriptstyle\rm Q$}\hbox{\raise
0.15\ht0\hbox to0pt{\kern0.4\wd0\vrule height0.7\ht0\hss}\box0}}
{\setbox0=\hbox{$\scriptscriptstyle\rm Q$}\hbox{\raise
0.15\ht0\hbox to0pt{\kern0.4\wd0\vrule height0.7\ht0\hss}\box0}}}}

\def\bbbz{{\mathchoice {\hbox{$\sf\textstyle Z\kern-0.4em Z$}}
{\hbox{$\sf\textstyle Z\kern-0.4em Z$}}
{\hbox{$\sf\scriptstyle Z\kern-0.3em Z$}}
{\hbox{$\sf\scriptscriptstyle Z\kern-0.2em Z$}}}}

\def\bbbc{{\mathchoice {\setbox0=\hbox{$\displaystyle\rm C$}\hbox{\hbox
to0pt{\kern0.4\wd0\vrule height0.9\ht0\hss}\box0}}
{\setbox0=\hbox{$\textstyle\rm C$}\hbox{\hbox
to0pt{\kern0.4\wd0\vrule height0.9\ht0\hss}\box0}}
{\setbox0=\hbox{$\scriptstyle\rm C$}\hbox{\hbox
to0pt{\kern0.4\wd0\vrule height0.9\ht0\hss}\box0}}
{\setbox0=\hbox{$\scriptscriptstyle\rm C$}\hbox{\hbox
to0pt{\kern0.4\wd0\vrule height0.9\ht0\hss}\box0}}}}

\begin{document}

\title{Homological index formulas for elliptic operators over $C^*$-algebras}

\author{Charlotte Wahl \footnote{This research was funded by a grant of AdvanceVT}}
\date{}

\maketitle

\begin{abstract}
We prove index formulas for elliptic operators acting between sections of $C^*$-vector bundles on a closed manifold. The formulas involve Karoubi's Chern character from $K$-theory of a $C^*$-algebra to de Rham homology of smooth subalgebras. We show how they apply to the higher index theorem for coverings and to flat foliated bundles, and prove an index theorem for $C^*$-dynamical systems associated to actions of compact Lie groups. In an Appendix we relate the pairing of odd $K$-theory and $KK$-theory to the noncommutative spectral flow and prove the regularity of elliptic pseudodifferential operators over $C^*$-algebras.
\end{abstract}

\section{Introduction}

One of the generalizations of the Atiyah-Singer index theorem is to elliptic pseudodifferential operators associated to $C^*$-vector bundles.  Mishenko--Fomenko introduced these operators and their index, an element in the $K$-theory of the $C^*$-algebra \cite{mf}. Furthermore they defined a Chern character for $C^*$-vector bundles and used it to formulate and prove an analogue of the Atiyah--Singer index theorem. However, in general it is not clear how to calculate the Mishenko--Fomenko Chern character of a $C^*$-vector bundle: Its definition is based on the map $K_0(C(M,\A))\ten \bbbc \to K_0(C(M))\ten  K_0(\A) \ten \bbbc \oplus K_1(C(M))\ten  K_1(\A) \ten \bbbc$ for a closed manifold $M$ and a unital $C^*$-algebra $\A$, which  exists by K\"unneth formula.

In this paper we prove index theorems for the same situation using Karoubi's Chern character from the $K$-theory of a $C^*$-algebra to the de Rham homology of smooth subalgebras \cite{ka}. 
Karoubi's Chern character is a generalization of the Chern character in differential geometry and is closely related to the Chern character in cyclic homology. Karoubi's de Rham homology has been used especially in noncommutative superconnections proof beginning with \cite{lo}.

We also prove (in the Appendix) that the pairing $K_1(\A) \times KK_1(\A,\B)\to K_1(\B)$, where $\A$, $\B$ are unital $C^*$-algebras, can be expressed in terms of the noncommutative spectral flow, which was introduced in the context of family index theory by Dai--Zhang \cite{dz}. See \cite{wa} for further references and a systematic account. The formula is well-known for $\B=\bbbc$ and the ordinary spectral flow. 

The main ingredient of the proof of the index theorem is a result about the compatibility of Karoubi's Chern character with the tensor product in $K$-theory. This allows the comparison of Karoubi's Chern character with  Mishenko-Fomenko's Chern character. 

Our proof generalizes the derivation of Atiyah's $L^2$-index theorem from the Mishenko--Fomenko index theorem in \cite{sc}. It is also closely related to the proof of an index theorem for flat foliated bundles in \cite{ji}, which is a special case of Connes' index theorem for foliated manifolds \cite[p. 273]{c} and implies the $C^*$-algebraic version of higher index theorem of Connes-Moscovici \cite{cm}. As an illustration we derive the the latter in detail from our formula. We also show how to apply the formula to flat foliated bundles. In this context we introduce a smooth subalgebra which is defined in more general situations than the one in \cite{ji}. 

We also prove an index theorem for Toeplitz operators associated to a $C^*$-dynamical system $(\A,G,\alpha)$ where $G$ is a compact Lie group. The Chern character involved here has been defined in \cite{co}. In \cite{le} a similar index theorem was proven for $G=\bbbr$ using Breuer-Fredholm operators. We relate both theorems in the case where the $\bbbr$-action is periodic.

In the Appendix we explain how the pairing of $K$-theory with $KK$-theory is related to  index theory and collect some useful facts about pseudodifferential operators over $C^*$-algebras beyond those proven in \cite{mf}, in particular that elliptic pseudodifferential operators  are adjointable as bounded operators between appropriate Sobolev spaces and regular as unbounded operators on a fixed Sobolev space. 

If not specified, tensor products between graded spaces are graded, and between Fr\'echet spaces they are completed projective.

{\it Acknowledgements:} I would like to thank Peter Haskell for helpful comments on previous versions of this paper.

\section{De Rham homology and the Chern character}

\subsection{Definition}

In this section we recall and slightly extend the definition of Karoubi's Chern character and collect properties that are relevant for index theory. The main reference is \cite{ka}.

Let $\Ai$ be a locally $m$-convex Fr\'echet algebra. 

The left $\Ai$-module of differential forms of order $k$ of $\Ai$ is defined as
$$\Ok \Ai:=\Ai \ten (\Ai/\bbbc)^{\ten k}$$ and  the $\bbbz$-graded space of all differential forms is $$\Oi\Ai :=\prod\limits_{k=0}^{\infty} \Ok \Ai \ .$$
There is a differential $\di$ on $\Oi\Ai$ of degree one defined by
$$\di(a_0 \ten \dots \ten a_k)=1 \ten a_0 \ten \dots \ten a_k$$ and a product determined by the properties that 
$$a_0 \ten \dots \ten a_k= a_0 \di a_1 \dots \di a_k$$ and that Leibniz rule holds which says that for $\alpha \in \Ok\Ai,~\beta \in \Oi\Ai$
$$\di(\alpha\beta)= (\di \alpha)\beta+ (-1)^k\alpha \di \beta \ .$$ 
With these structures $\Oi\Ai$ is a graded differential locally $m$-convex Fr\'echet algebra.
 
For a closed manifold $M$
$$\Ob^{p,q}(M,\Ai):=\Ob^p(M,\Ob_q\Ai)=\Ob^p(M) \ten \Ob_q\Ai \ ,$$ where $\Ob^*(M)$ is the space of smooth differential forms on $M$. 

We call an open subset $U \subset M$ regular if the compactly supported de Rham cohomology $H_c^*(U)$ is finite-dimensional and if there are open subsets $U_0,U_1$ with $\ov{U_0} \subset U$ and $\ov{U} \subset U_1$ such that there is a smooth homotopy $F:[0,1] \times U_1 \to U_1$ with $F(0,x)=x$ for all $x \in U$ and such that $F^{-1}(\{1\} \times U) \subset \{1\} \times U_0$ and $F^{-1}(\{t\} \times U) \subset \{t\} \times U$ for all $t \in [0,1]$.  

For a regular open subset $U$ in $M$ we define $\Ob^{p,q}_0(U,\Ai)$ to be the closure of the subspace of $\Ob^{p,q}(M,\Ai)$ spanned by forms with support in $U$. 

The product on $\Ob^{**}_0(U,\Ai)$ is determined by the natural isomorphism $\Ob_0^{**}(U,\Ai) \cong \Ob^*_0(U) \ten \Ob_*\Ai$. Here the right hand side is understood as a graded tensor product of graded algebras. Let $d_U$ be the de Rham differential on $U$. The differential of the total complex of the double complex $(\Ob_0^{**}(U,\Ai), d_U,\di)$ is denoted by $d_{tot}$ and its homology by   
$H_0^*(U,\Ai)$. The definition does not depend on the embedding of $U$ into $M$ as a regular subset.

For a closed manifold $M$ we usually omit the suffix and write $H^*(M,\Ai)$.

For $F$ as above let $f_t=F(t,\cdot): U \to U$. Then $f_1^*:H_0^*(U) \to H_c^*(U)$ is inverse to the map $H_c^*(U) \to H_0^*(U)$, since for a closed form $\omega \in \Ob_0^*(U)$ the form $f_1^*\omega$ is a closed form supported in $U$ and $f_1^*\omega-\omega = d_U \int_0^1 F^*\omega$. Hence $H_c^*(U) \cong H_0^*(U)$.

The isomorphism $\Ob_0^{**}(U,\Ai) \cong \Ob_0^*(U) \ten \Ob_*\Ai$ induces isomorphisms 
$$\Ob_0^{**}(U,\Ai)/\ov{[\Ob_0^{**}(U,\Ai),\Ob_0^{**}(U,\Ai)]_s}\cong\Ob_0^*(U)\ten \Oi\Ai/\ov{[\Oi\Ai,\Oi\Ai]_s} \ ,$$ 
$$H_0^n(U,\Ai)\cong \oplus_{p+q=n} H_0^p(U) \ten H_0^q(\Ai) \cong \oplus_{p+q=n}H_0^p(U,H_q(\Ai)) \ .$$
These isomorphisms have been proven in \cite[\S\S 4.7, 4.8]{ka} in a slighly different situation. The proof carries over. It uses completed tensor products, therefore we use $\Ob_0^*(U)$ instead of compactly supported forms for the definition of cohomology.  The proof uses furthermore the fact that $H_0^*(U)$ is finite-dimensional.
 
We call a smooth possibly noncompact manifold $M$ regular if there is a covering $(U_n)_{n \in \bbbn}$ by regular subsets with $U_n \subset U_{n+1}$. 

Extending a form by zero induces a well-defined push forward map $H_0^*(U_n,\Ai) \to H_0^*(U_{n+1},\Ai)$ so that we can define
$$H_c^*(M,\Ai)=\varinjlim\limits_{n \to \infty} H_0^*(U_n,\Ai) \ .$$

It is clear that $H_c^*(M)$ agrees with the compactly supported de Rham cohomology of $M$.

If $M \to B$ is a fiber bundle of regular oriented manifolds, then integration over the fiber yields a homomorphism
$$\int_{M_b}:H_c^{*}(M,\Ai) \to H_c^{*-\dim M_b}(B,\Ai)  . $$

The de Rham homology of $\Ai$ is $$H_*(\Ai):=H^*(*,\Ai)\ ,$$ where $*$ is the point.

If $M$ is a closed manifold, we usually write $H^*(M,\Ai)$ for $H^*_0(M,\Ai)$. 
Note that then the quotient map $\Ob_n(\C(M,\Ai)) \to \oplus_{p+q=n}\Ob^{p,q}(M,\Ai)$ induces a homomorphism 
$$H_*(\C(M,\Ai)) \to H^{*}(M,\Ai) \ .$$

We proceed with the definition and the properties of the Chern character.

Let $\A$ be a unital $C^*$-algebra and let $\Ai \subset \A$ be a dense subalgebra that is closed under involution and holomorphic functional calculus in $\A$. Assume that $\Ai$ is endowed with the topology of a locally $m$-convex algebra such that $\Ai \incl \A$ is continuous. We call such a subalgebra  a smooth subalgebra of $\A$. 

Let $M$ be a regular manifold. Recall that $K_0(C_0(M,\A))=\Ker (K_0(C_0(M,\A)^+) \to K_0(\bbbc))$, where $C_0(M,\A)^+$ denotes the unitalization of $C_0(M,\A)$. Since $\C_c(M,\Ai)^+$ is dense and closed under holomorphic functional calculus in $C_0(M,\A)^+$, we have that $K_0(\C_c(M,\Ai)) \cong K_0(C_0(M,\A))$.

The Chern character form of a projection $P \in M_n(\C_c(M,\Ai)^+)$ is defined as 
$$\ch^M_{\Ai}(P):=\sum_{k=0}^{\infty}\frac{(-1)^k}{(2\pi i)^k k!}\tr P(d_{tot} P)^{2k} \ .$$ The normalization differs from the normalization in \cite{ka} and is chosen such that the Chern character of the Bott element $B \in K_0(C_0((0,1)^2))$ integrated over $(0,1)^2$ equals $1$. (There is also some ambiguity about the sign of the Bott element $B$ in the literature. Here we take $B=1-[H] \in \Ker(K_0(C(S^2)) \to K_0(\bbbc))$, where $H$ is the Hopf bundle.)

In the following proposition we denote by $P_{\infty} \in M_n(\bbbc)$ the image of $P \in M_n(\C_c(M,\Ai)^+)$ under ``evaluation at infinity''.

\begin{prop}
\begin{enumerate}
\item $\ch^M_{\Ai}(P)$ is closed.
\item
Let $P:[0,1] \to M_n(\C_c(M,\Ai)^+)$ be a differentiable path of projections and let $U \subset M$ be such that $\supp(P(t)-P_{\infty}(t)) \subset U$ for all $t \in [0,1]$. Then there is a form $\alpha \in \Ob^{**}_0(U,\Ai)/\ov{[\Ob^{**}_0(U,\Ai),\Ob_0^{**}(U,\Ai)]_s} $ such that $d_{tot} \alpha = \ch^M_{\Ai}(P(1))- \ch^M_{\Ai}(P(0))$.
\item The Chern character form induces a well-defined homomorphism $$K_0(C_0(M,\A)) \to H_c^*(M,\Ai) \ .$$
\end{enumerate}
\end{prop}

\begin{proof}
For $M$ compact the proofs are standard. We include the proof of 2) in order to show that it works in the noncompact case as well:

From Leipniz rule on deduces that the terms $PP'P,~(1-P)P'(1-P),~P(d_{tot} P)P,~(1-P)(d_{tot} P)(1-P)$ all vanish. 

Hence
\begin{eqnarray*}
 \tr (P(d_{tot}P)^{2k})'
 &=& \tr P'(d_{tot}P)^{2k}+  \tr P((d_{tot}P)^{2k})' \\
&=& \tr P((d_{tot}P)^{2k})'\\
&=& \sum\limits_{i=0}^{2k-1} \tr P(d_{tot} P)^i(d_{tot} P)'(d_{tot} P)^{2k-i-1} \ .
\end{eqnarray*}

This vanishes for $k=0$.

For $i$ even and $k \neq 0$
\begin{eqnarray*}
\lefteqn{\tr P(d_{tot} P)^i(d_{tot} P)'(d_{tot} P)^{2k-i-1}}\\
 &=& \tr (d_{tot} P)^i P(d_{tot} P')(d_{tot} P)^{2k-i-1}\\
&=& \tr (d_{tot} P)^i (d_{tot} (PP'))(d_{tot} P)^{2k-i-1}- \tr (d_{tot} P)^i (d P)P'(d_{tot} P)^{2k-i-1}\\
&=&\tr (d_{tot} P)^i (d_{tot} (PP'))(d_{tot} P)^{2k-i-1}\\
&=&d_{tot} \tr P(d_{tot} P)^{i-1}(d_{tot} (PP'))(d_{tot} P)^{2k-i-1} \ .
\end{eqnarray*}
 Note that $\tr P(d_{tot} P)^{i-1}(d_{tot} (PP'))(d_{tot} P)^{2k-i-1}$ vanishes on $U$ for $k \neq 0$. 

For $i$ odd the argument is similar.
\end{proof}

We define the odd Chern character via the following diagram:
\begin{eqnarray} \label{bott1}  \begin{CD} 
K_0(C_0((0,1) \times M, \A)) @>  \cong >> K_1(C_0(M,\A)) \\
@VV   \ch^{(0,1) \times M}_{\Ai}   V@VV \ch^M_{\Ai} V\\
H^{ev}_c((0,1) \times M,\Ai) @> \int_0^1 >> H^{odd}_c(M,\Ai) \ .
\end{CD}  
\end{eqnarray}

Note that  $\int_0^1:H_c^*((0,1) \times M, \Ai) \to H_c^*(M,\Ai)$ is an isomorphism by $H_c^*((0,1)\times M, \Ai) \cong H_c^*((0,1)) \ten H_c^*(M,\Ai)\cong H_c^*(M,\Ai) \ .$

In the following we derive  a formula for the odd Chern character. It is analogous to those well-known in de Rham cohomology and cyclic homology (compare with \cite{gl}).

\begin{prop} For $u \in U_n(\C_c(M,\Ai)^+)$ in $H^*_c(M,\Ai)$
$$\ch^M_{\Ai}([u])= \sum_{k=1}^{\infty}\left(\frac{-1}{2\pi i}\right)^k\frac{(k-1)!}{(2k-1)!} u^*(d_{tot} u)((d_{tot} u^*)(d_{tot} u))^{k-1}  \ .$$
\end{prop}

\begin{proof}
Here we use that the Chern character can be defined in terms of noncommutative connections \cite{ka}.

Let $P_n \in M_{2n}(\bbbc)$ be the projection onto the first $n$ components. Let $W(t) \in \C([0,1],U_{2n}(\C_c(M,\Ai)^+))$ with $W(0)=1$ and  $W(1)=\diag(u,u^*)$. Then the isomorphism $K_1(C_0(M,\A)) \to K_0(C((0,1) \times M,\A))$ maps $[u]$ to $[WP_nW^*]-[P_n]$. The Chern character is independent of the choice of the connection \cite[Th. 1.22]{ka}, thus we may use the connection $WP_n(d_{tot} + dx~\ra_x + xW(1)^*d_{tot}(W(1)))W^*$ on the projective $\C_c((0,1) \times M,\Ai)^+$-module $WP_nW^*(\C_c((0,1) \times M,\Ai)^+)^n$ for its calculation. It follows that 
\begin{eqnarray*}
\lefteqn{\ch_{\Ai}^M(u)}\\
&=&\ch_{\Ai}^{(0,1)\times M}(WP_nW^*)\\
&=& \sum_{k=0}^{\infty}\int_0^1 \frac{(-1)^k}{(2 \pi i)^k k!}\tr(x^2 u^*(d_{tot} u)u^*(d_{tot} u)+ x (d_{tot} u^*)(d_{tot} u) +dx ~u^*(d_{tot} u))^k \\
&=&  \sum_{k=0}^{\infty}\frac{(-1)^k}{(2\pi i)^k k!} \int_0^1 \tr((x-x^2)(d_{tot} u^*)(d_{tot} u) +dx ~u^*(d_{tot} u))^k \\
&=&\sum_{k=1}^{\infty}\frac{(-1)^k}{(2 \pi i)^k k!} \int_0^1 dx~(x-x^2)^{k-1}u^*(d_{tot} u)((d_{tot} u^*)(d_{tot} u))^{k-1}\\
&=&\sum_{k=1}^{\infty}\left(\frac{-1}{2\pi i}\right)^k\frac{(k-1)!}{(2k-1)!} u^*(d_{tot} u)((d_{tot} u^*)(d_{tot} u))^{k-1} \ .
\end{eqnarray*}
\end{proof}

\subsection{Chern character and tensor products}

From now on assume that $M$ is a closed manifold.

Let $K_i(\A)_{\bbbc}:= K_i(\A) \ten \bbbc$. 

In the following we prove the compatibility of the Chern character with the Bott periodicity map $K_1(C_0((0,1),\A) \cong K_0(\A)$ and with the K\"unneth formulas $$K_0(C(M))_{\bbbc} \ten K_0(\A)_{\bbbc} \oplus K_1(C(M))_{\bbbc} \ten K_1(\A)_{\bbbc}\cong K_0(C(M,\A))_{\bbbc}$$
and
$$K_0(C(M))_{\bbbc} \ten K_1(\A)_{\bbbc} \oplus K_1(C(M))_{\bbbc} \ten K_0(\A)_{\bbbc}\cong K_0(C(M,\A))_{\bbbc} \ .$$
These isomorphisms are defined via the tensor product $$K_i(C(M)) \ten K_j(\A) \to K_{i+j}(C(M,\A)),~i,j \in \bbbz/2 \ .$$
The tensor product is injective, hence we may consider $K_i(C(M)) \ten K_j(\A)$ as a subspace of $K_{i+j}(C(M,\A))$.

First recall the definition of the tensor product.
For $i,j=0$ the tensor product is induced by the tensor product of projections. The remaining three cases are derived from the tensor product of projections using Bott periodicity, for example
\begin{eqnarray*}
K_1(C(M)) \ten K_1(\A) &\cong& K_0(C_0((0,1) \times M)) \ten K_0(C_0((0,1))\ten \A)\\
& \to & K_0(C_0((0,1)^2 \times M,\A)) \\
& \cong& K_0(C(M,\A)) 
\end{eqnarray*}
and
\begin{eqnarray*}
K_0(C(M)) \ten K_1(\A) &\cong& K_0(C(M)) \ten K_0(C_0((0,1), \A))\\
& \to & K_0(C_0((0,1) \times M,\A)) \\
& \cong& K_1(C(M,\A)) \ .
\end{eqnarray*}

A standard calculation (see \cite[Th. 1.26]{ka}) shows that the tensor product for $i=j=0$ is compatible with the Chern character, namely for $a \in K_0(C(M))$ and $b \in K_0(\A)$  $$\ch^M_{\Ai}(a \ten b)=\ch^M(a)\ch_{\Ai}(b) \ .$$ 
In the following proposition  $\beta:K_0(C(M,\A)) \to K_0(C_0((0,1)^2 \times M,\A),~ a \mapsto a \ten B$ is the Bott periodicity map.  

\begin{prop} 
\begin{enumerate}
\item For
$a \in K_0(C_0((0,1)^2 \times M, \A))$
$$\ch^M_{\Ai} \beta^{-1}(a)= \int_{(0,1)^2} \ch^{(0,1)^2 \times M}_{\Ai} (a) \ .$$
\item For $a \in K_0(C_0((0,1) \times M))$ and $b \in K_0(C_0((0,1),\A))$ 
$$\ch_{\Ai}^M \beta^{-1}(a\ten b)=\int_{(0,1)^2} \ch^{(0,1)\times M}(a) \ch^{(0,1)}_{\Ai} (b) \ .$$
\end{enumerate}
\end{prop}

\begin{proof} We consider $K_0(C_0((0,1)^2\times M,\A))$ as a subgroup of $K_0(C(T^2 \times M,\A))$.

1) Let $b \in K_0(C(M,\A))$ with $a=B \ten b$.  Then 
\begin{eqnarray*}
\int_{(0,1)^{2}} \ch^{(0,1)^2\times M}_{\Ai}(B \ten b) &=& \int_{T^{2}} \ch^{T^{2} \times M}_{\Ai}(B \ten b) \\
&=&\ch^M_{\Ai}(b) \int_{T^{2}}\ch^{T^{2}}(B) \\
&=& \ch^M_{\Ai}(b) \ .
\end{eqnarray*}

2) The assertion follows from the commutative diagram 
$$\begin{CD}
K_0(C_0((0,1)\times M)) \ten K_0(C_0((0,1), \A)) @>   >> K_0(C(S^1 \times M)) \ten K_0(C(S^1,\A))\\
@VV      V@VV \ten V\\
K_0(C_0((0,1)^2\times M, \A)) @> \incl >> K_0(C(T^2 \times M,\A)) \\
@VV \ch^M_{\Ai} \circ \beta^{-1} V@VV \int_{T^2}\ch^{T^2\times M}_{\Ai} V\\
H^*(M,\Ai) @> = >> H^*(M,\Ai)
\end{CD} \ .$$
Since the horizontal arrows are inclusions, the first vertical map on the left hand side is determined by the first vertical map on the right hand side. The second square commutes by 1). 
\end{proof}

\begin{cor}
The diagram 
$$\begin{CD}
K_1(C_0((0,1) \times M, \A) @> \cong  >> K_0(C (M,\A)) \\
@VV   \ch^{(0,1) \times M}_{\Ai} V@VV  \ch^M_{\Ai} V\\
H^{odd}_c((0,1) \times M,\Ai) @> \int_0^1  >> H^{ev}(M,\Ai) 
\end{CD}$$
commutes. 
\end{cor}

\begin{proof}
Consider the diagram
$$\begin{CD}
K_1(C_0((0,1) \times M, \A)) @< \cong  << K_0(C_0((0,1)^2 \times M,\A)) @> \beta^{-1} >> K_0(C(M,\A)) \\
@VV   \ch^{(0,1)\times M}_{\Ai} V@VV  \ch^{(0,1)^2\times M}_{\Ai}  V@VV \ch^M_{\Ai} V\\
H_c^{odd}((0,1)\times M,\Ai) @< \int_0^1 << H^{ev}_c((0,1)^2\times M,\Ai)  @> \int_{(0,1)^2}  >> H^{ev}(M,\Ai) 
\end{CD} \ .$$
The first square commutes by diagram \ref{bott1} applied to $(0,1) \times M$. The second square commutes by the first part of the previous proposition.
\end{proof}

We denote by $$R_{jk}:K_i(C(M,\A))_{\bbbc} \to K_j(C(M))_{\bbbc}\ten K_k(\A)_{\bbbc}\subset K_i(C(M,\A))_{\bbbc}$$ the projections induced by the K\"unneth formulas.

We have a tensor product $$\ch^M \ten \ch_{\Ai}:K_i(C(M))_{\bbbc} \ten K_i(\A)_{\bbbc} \to H^*(M) \ten H_*(\Ai)\cong H^*(M,\Ai) \ .$$

\begin{prop} 
\label{compten}
\begin{enumerate}
\item On $K_0(C(M,\A))_{\bbbc}$ 
$$\ch^M_{\Ai}= (\ch^M \ten \ch_{\Ai}) \circ R_{00} + (\ch^M \ten \ch_{\Ai}) \circ R_{11} \ .$$
 
\item On $K_1(C(M,\A))_{\bbbc}$ 
$$\ch^M_{\Ai}=(\ch^M \ten \ch_{\Ai}) \circ R_{01} + (\ch^M \ten \ch_{\Ai}) \circ R_{10} \ .$$ 
\end{enumerate}
\end{prop}

\begin{proof} 1) follows from the previous proposition:
Let $a \ten b \in K_0(C(M,\A))$ with $a \in K_1(C(M)))$ and $b \in K_1(\A)$. Let
$a$ correspond to $\tilde a \in K_0(C_0((0,1) \times M))$ and $b$ to $\tilde b \in K_0(C_0((0,1),\A))$.

Then by definition $\ch^M(a)=\int_0^1\ch^{(0,1)\times M}(\tilde a)$ and  $\ch_{\Ai}(b)=\int_0^1\ch^{(0,1)}_{\Ai}(\tilde b)$. Now by the previous lemma 
\begin{eqnarray*}
\ch^M_{\Ai}(a \ten b)&=& \int_{(0,1)^2} \ch^{(0,1)\times M}(\tilde a) \ch^{(0,1)}_{\Ai} (\tilde b) \\
&=&\int_0^1\ch^{(0,1)\times M}(\tilde a) \int_0^1 \ch^{(0,1)}_{\Ai} (\tilde b)\\
&=& \ch^M(a) \ch_{\Ai}(b) \ .
\end{eqnarray*}

2) follows applying by 1) to $K_0(C_0((0,1) \times M,\A))$ since the Chern character interchanges the suspension isomorphisms in $K$-theory and de Rham homology.
\end{proof}

Define $\Ch_{\A}^M$ as the map
\begin{eqnarray*}
K_0(C(M,\A)) &\to& K_0(C(M,\A))_{\bbbc} \\
&\cong& K_0(C(M))_{\bbbc} \ten K_0(\A)_{\bbbc} \oplus K_1(C(M))_{\bbbc} \ten K_1(\A)_{\bbbc}\\
&\stackrel{\ch^M}{\lr}& H_{dR}^{ev}(M) \ten K_0(\A) \oplus H_{dR}^{odd}(M) \ten K_1(\A) \ .
\end{eqnarray*}
and analogously for $K_1(C(M,\A))$.
This is the Chern character introduced by Mishenko-Fomenko \cite{mf}. 

The previous Proposition is equivalent to the equation 
\begin{eqnarray} \ch_{\Ai} \circ \Ch_{\A}^M&=& \ch_{\Ai}^M  \ . \label{ChMF} \end{eqnarray}

\subsection{Pairing with cyclic cocycles}

In the noncommutative geometry the Chern character with values in the cyclic homology is more common than the one with values in the de Rham homology. De Rham homology can be paired with normalized cyclic cocycles; in this pairing both Chern characters agree up to normalization:

Let $\ov{C}_n^{\lambda}(\Ai)$ be the quotient of the algebraic tensor product $(\Ai/\bbbc)^{\ten n+1}$ by the action of $\bbbz/(n+1)\bbbz$.
Let $$b:\ov{C}_n^{\lambda}(\Ai) \to \ov{C}_{n-1}^{\lambda}(\Ai) \ ,$$
\begin{eqnarray*}
b(a_0 \ten \dots a_n)&=&(-1)^na_na_0 \ten \dots \ten a_{n-1}\\
&& + \sum_{i=0}^{\infty}(-1)^i
a_0 \ten \dots  \ten a_ia_{i+1}\ten \dots a_n \ .
\end{eqnarray*} 
The homology of the complex $(\ov{C}_*^{\lambda}(\Ai),b)$ is the reduced  cyclic homology $\ov{HC}_*(\Ai)$. Using the completed projective tensor product instead of the algebraic one we obtain the topological reduced cyclic homology $\ov{HC}_*^{top}(\Ai)$.  Furthermore we denote by $\ov{HC}_*^{sep}(\Ai)$ the topological homology of $(\ov{C}_n^{\lambda}(\Ai),b)$, i.e. we use the completed projective tensor product and quotient out the closure of the range of $b$. 

The reduced cyclic cohomology $\ov{HC}^*(\Ai)$ is the homology of the dual complex $(\ov{C}_{\lambda}^n(\Ai),b^t)$ (in the algebraic sense). Elements of $\ov{C}_{\lambda}^n(\Ai)$ are called normalized cochains. The continuous reduced cyclic cohomology $\ov{HC}_{top}^*(\Ai)$ is the homology of the topological dual complex.

The pairing  $\ov{HC}_{top}^*(\Ai) \ten \ov{HC}_*^{top}(\Ai) \to \bbbc$ descends to a pairing $\ov{HC}_{top}^*(\Ai) \ten \ov{HC}_*^{sep}(\Ai) \to \bbbc$. Furthermore the quotient map $\Ob_n(\Ai) \to \ov{C}^n_{\lambda}(\Ai)$ induces an homomorphism $H_n(\Ai) \to \ov{HC}_n^{sep}(\Ai)$, which is an embedding for $n \ge 1$ (see \cite[\S\S 4.1 and 2.13]{ka}). In degree zero there is a pairing of $H_0(\Ai)=\Ai/[\Ai,\Ai]$ with traces on $\Ai$. 

The Chern character $\ch^{\lambda}:K_0(\Ai) \to \ov{HC}_*(\Ai)$ is defined by 
$$\ch^{\lambda}(p)= \sum_{m=0}^{\infty}(-1)^m \tr p^{\ten 2m+1}$$ for a projection $p \in M_n(\Ai)$. Hence the composition
$$K_0(\Ai) \stackrel{\ch^{\lambda}}{\lr} \ov{HC}_*(\Ai) \to \ov{HC}_*^{sep}(\Ai)$$ agrees up to normalization with the map
$$K_0(\Ai) \stackrel{\ch_{\Ai}}{\lr} H_*(\Ai) \incl \ov{HC}_*^{sep}(\Ai) \ .$$ In particular if $\phi \in \ov{HC}_{top}^m(\Ai)$, then  $$\phi\circ \ch^{\lambda}=(2 \pi i)^m m!~ \phi \circ \ch_{\Ai} \ .$$

\section{Index theorems}

In the following we give a formulation
of the Mishenko--Fomenko index theorem, which is different from the original one and adapted to the applications. Furthermore we translate its proof in the language of $KK$-theory: We show the compatibility of the Chern character with the pairing $K_i(C(M,\A)) \ten KK_j(C(M),\bbbc) \to K_{i+j}(\A)$ for $i,j \in \bbbz/2$, where on $KK_j((C(M),\bbbc)$ we use the Chern character from $K$-homology to de Rham homology of $M$. We refer to Appendix \ref{indKK} for some facts about the connection of $KK$-theory to index theory.

\begin{lem}\begin{enumerate}
\item For $x \in K_i(C(M)) \ten K_i(\A) \subset K_0(C(M,\A))$ and $y \in KK_j(C(M),\bbbc)$ with $i \neq j$ 
$$x \ten_{C(M)} y=0 \in K_j(\A) \ .$$
\item For $x \in K_i(C(M)) \ten K_j(\A) \subset K_1(C(M,\A))$ and $y \in KK_j(C(M),\bbbc)$ with $i \neq j$ $$x \ten_{C(M)} y=0 \in K_i(\A) \ .$$
\end{enumerate}
It follows that for $x \in K_i(C(M,\A))$ and $y \in KK_j(C(M),\bbbc)$
$$x\ten_{C(M)} y=R_{j,i+j}(x)\ten_{C(M)} y \in K_{i+j}(\A)_{\bbbc} \ .$$
\end{lem}

\begin{proof}
1) Let $B_1 \in KK_0(C_0((0,1)^2),\bbbc)$ be the Bott element. By the standard isomorphism $K_i(\A) \cong KK_i(\bbbc,\A)$ and the fact that the tensor product in $K$-theory is a special case of the Kasparov product all we have to show is that for $a \in KK_j(\bbbc,C_0((0,1),\A))$ and $b\in KK_j(\bbbc,C_0((0,1) \times M))$
$$((a \ten b)\ten_{C_0((0,1)^2)} B_1) \ten_{C(M)} y=0 \ .$$
This follows from the associativity of the product and the fact that $(b\ten_{C_0((0,1))} B_1) \ten_{C(M)} y \in KK_0(C_0((0,1)),\bbbc)=0$.

2) Let $i=0,~j=1$. Let $B_2 \in KK_1(C_0((0,1)),\bbbc)$ be the Bott element.  Let $a  \in K_0(C_0((0,1),\A))$ and $b\in K_0(C(M))$. 
Then $((a \ten b) \ten_{C_0((0,1)^2)} B_2)\ten_{C(M)}  y=0$ by associativity and since $(b \ten_{C_0((0,1))} B_2) \ten_{C(M)}  y \in KK_0(C_0((0,1)),\bbbc)=0$.
The proof for $i=1,~j=0$ is analogous.
\end{proof}

Let now $\chh_M:KK_i(C(M),\bbbc) \to H_*(M)$ be the homological Chern character where $H_*(M)$ is the de Rham homology of $M$ with complex coefficients.

For the following proposition note that the pairing $\langle~,~\rangle: H^*(M) \times H_*(M) \to \bbbc$ induced a pairing $\langle~,~\rangle:(H^*(M) \ten K_*(\A)) \times H_*(M) \to K_*(\A)_{\bbbc}$.

\begin{lem}
For $x \in K_i(C(M,\A))$ and $y \in KK_j(C(M),\bbbc)$
$$\langle \Ch_{\A}^M x, \chh_M y\rangle =\langle \Ch_{\A}^M R_{j,i+j}(x), \chh_M y \rangle \in K_{i+j}(\A)_{\bbbc} \ .$$
\end{lem}

\begin{proof} Consider the case $i,j=0$. We have to show that $\langle \Ch_{\A}^M R_{11}(x), \chh_M y \rangle =0$ or equivalently that $\langle \Ch_{\A}^M (x), \chh_M y \rangle =0$ for $x=x_1 \ten x_2$ with $x_1 \in K_1(C(M))$ and $x_2 \in K_1(\A)$. 

Clearly $\Ch^M_{\A} x=(\ch^M x_1) x_2$. 

Hence 
$$\langle \Ch_{\A}^M(x), \chh_M(y) \rangle= \langle \ch^M x_1, \chh_M(y)\rangle x_2 \ .$$
Since $\ch^Mx_1 \in H^{odd}(M)$ and $\chh_M(y) \in H_{ev}(M)$, the pairing $\langle \ch^M x_1, \chh_M(y)\rangle$ vanishes. 

The remaining three cases are analogous.
\end{proof}

\begin{prop}
If $x \in K_i(C(M,\A))$ and $y \in KK_j(C(M),\bbbc)$, then 
  $$x\ten_{C(M)} y = \langle \Ch^M_{\A}x , \chh_M y\rangle \in K_{i+j}(\A)_{\bbbc} \ .$$  
\end{prop}

\begin{proof}
By the first lemma $x \ten_{C(M)}y=R_{j,i+j}(x) \ten_{C(M)}y \in K_{i+j}(\A)_{\bbbc}$. By the previous lemma the right hand side of the formula also only depends on $R_{j,i+j}(x)$. Therefore and by linearity we may restrict to the case where $x=x_1 \ten x_2$ with $x_1\in K_j(C(M))$ and $x_2 \in K_{i+j}(\A)$.
Then in $K_{i+j}(\A)_{\bbbc}$
\begin{eqnarray*}
x \ten_{C(M)}y &=& (x_1\ten_{C(M)}y) x_2  \\
&=& \langle \ch^M x_1, \chh_M y\rangle x_2 \\
&=&\langle \Ch^M_{\A}x, \chh_M y \rangle \ .
\end{eqnarray*} 
\end{proof}

Using formula \ref{ChMF} and considering the pairing $\langle ~,~\rangle: H^*(M,\Ai) \times H_*(M) \to H_*(\Ai)$ we obtain:

\begin{cor}
If $x \in K_i(C(M,\A))$ and $y \in KK_j(C(M),\bbbc)$, then 
  $$\ch_{\Ai}(x\ten_{C(M)} y) = \langle \ch^M_{\Ai}x , \chh_M y\rangle \in H_*(\Ai) \ .$$  
\end{cor}

In the following we translate these results into a more classical language (see Appendix \ref{indKK}):

Now let $M$ be a closed Riemannian manifold and let $E$ be a  hermitian, possibly $\bbbz/2$-graded, complex vector bundle on $M$.

Let $D:\C(M,E)\to \C(M,E)$ be an elliptic symmetric pseudodifferential operator of order $1$. If $E$ is graded, then $D$ is assumed to be odd.
In the ungraded case the symbol $\sigma(D)$ defines an element in $K_1(C_0(T^*M))$, in the graded case $[\sigma(D^+)]\in K_0(C_0(T^*M))$. If $E$ is ungraded, then $[(L^2(M,E),D)] \in KK_1(C(M),\bbbc)$, else $[(L^2(M,E),D)] \in KK_0(C(M),\bbbc)$. In the ungraded case the values of the index $\ind$ are in $K_1(\A)$, in the graded case in $K_0(\A)$.

Define the $\A$-vector bundle $L(U):=([0,1] \times M\times \A^n)/(0,x,v)\sim (1,x,U(x)v)$ on $S^1 \times M$ and let $\dira_{L(U)}$ be the  operator $\frac{1}{i} \frac{d}{dx}$ acting on the sections of $L(U)$. Pull $E$ back to $S^1 \ten M$. Then $\phi(t)D+(1-\phi(t))UDU^*$ is well-defined on $L^2(S^1 \times M,L(U) \ten E)$. 

Let $\pi_!:H^*_c(TM) \to H^*(M)$ be integration over the fiber and $k=\frac{\dim M(\dim M+1)}{2}$. 
 
\begin{theorem}
\label{indtheor}
\begin{enumerate}
\item Let $P\in M_n(\C(M,\A))$ be a projection. 
\begin{enumerate}
\item Assume that $E$ is $\bbbz/2$-graded. Then
$$\ch_{\Ai} \ind P(\oplus^n D^+)P = (-1)^k\int_M \Td(M) \pi_!\ch^{TM} [\sigma(D^+)] \ch^M_{\Ai}[P] \ .$$
\item If $E$ is ungraded, then
$$\ch_{\Ai}\ind P(\oplus^n D)P = (-1)^k\int_M \Td(M) \pi_!\ch^{TM} [\sigma(D)] \ch^M_{\Ai}[P]   \ .$$
\end{enumerate}
\item Let $U \in U_n(\C(M,\A))$ be a unitary.
\begin{enumerate}
\item If $E$ is ungraded, then 
$$\ch_{\Ai}\spfl((1-t)D +tUDU^*)=(-1)^k\int_M \Td(M) \pi_!\ch^{TM}[\sigma(D)] \ch^M_{\Ai}[U]   \ .$$
\item If $E$ is $\bbbz/2$-graded and $\sigma$ is the grading operator, then
$$\ch_{\Ai}\ind(-\sigma\dira_{L(U)}+ i\sigma (\chi(t)D+(1-\chi(t))UDU^*))$$
$$=(-1)^k\int_M \Td(M) \pi_!\ch^{TM}[\sigma(D^+)] \ch^M_{\Ai}[U]  \ .$$
\end{enumerate}
\end{enumerate}
\end{theorem}

See Appendix \ref{indKK} for more possibilities to express the left hand side of 2(a) and 2(b).

\section{Applications}

\subsection{Higher index theory for coverings and flat foliated bundles}

In the following we deduce the higher index theorem for coverings of Connes-Moscovici \cite{cm} from the previous formulas. 
We do not recover the theorem in full generality (which calculates the pairing of an index in algebraic $K$-theory with group cocycles), but for extendable cocycles. 

Let $\Gamma$ be a discrete group.
 
We begin by recalling some facts about the group cohomology $H^*(\Gamma)$, in particular how to embed it into $\ov{HC}^*(\Gamma)$. 

Let $$C^n(\Gamma)=\{\tau:\Gamma^{n+1} \to \bbbc,~\tau(gg_0,\dots,gg_n)=\tau(g_0,\dots,g_n) \mbox{ for all } g \in \Gamma\} $$ and let $$d_{\Gamma}:C^n(\Gamma)\to C^{n+1}(\Gamma),$$
$$d_{\Gamma}\tau(g_0,\dots, g_{n+1})=\sum_{j=0}^{n+1} (-1)^j\tau(g_0, \dots,g_{j-1},g_{j+1}, \dots g_{n+1})\ .$$ 
We denote by $C_{al}^n(\Gamma)\subset C^n(\Gamma)$ the subspace of alternating elements. The homology of $(C_{al}^*(\Gamma),d_{\Gamma})$ is $H^*(\Gamma)$ (as is the homology of $(C^*(\Gamma),d_{\Gamma})$).

Furthermore let $$\ov{C}_{\lambda}^n(\bbbc\Gamma)_{<e>}=\{c \in \ov{C}_{\lambda}^n(\bbbc\Gamma)~|~ c(g_0,g_1, \dots g_n)=0 \mbox{ for } g_0g_1 \dots g_n \neq e \}$$  and let  $$\ov{C}_{\lambda}^n(\bbbc\Gamma)_{<g>\neq <e>} = \{c \in \ov{C}_{\lambda}^n(\bbbc\Gamma)~|~ c(g_0,g_1, \dots g_n)=0 \mbox{ for } g_0g_1 \dots g_n = e \} \ .$$
The complex $(\ov{C}_{\lambda}^*(\bbbc\Gamma),b^t)$ decomposes into a direct sum $(\ov{C}_{\lambda}^*(\bbbc\Gamma)_{<e>},b^t) \oplus (\ov{C}_{\lambda}^*(\bbbc\Gamma)_{<g>\neq <e>},b^t)$. 

In the following we assume $n \ge 1$.

For $c \in \ov{C}_{\lambda}^n(\bbbc\Gamma)_{<e>}$ define $\tau_c \in C_{al}^n(\Gamma)$ by $$\tau_c(e,g_1, \dots, g_n) := c(g_n^{-1}, g_1,g_1^{-1}g_2,g_2^{-1}g_3,\dots g_{n-1}^{-1}g_n)$$
and for $\tau \in C^n_{al}(\Gamma)$ define $c_{\tau} \in \ov{C}_{\lambda}^n(\bbbc\Gamma)_{<e>}$ by $$c_{\tau}(g_0,g_1, \dots,g_n)=\tau(e,g_1,g_1g_2, \dots, g_1\dots g_n) \mbox{ if } g_0g_1 \dots g_n=e$$ and $$c_{\tau}(g_0,g_1, \dots,g_n)=0 \mbox{ if } g_0g_1 \dots g_n \neq e \ .$$
The maps 
$$\ov{C}_{\lambda}^n(\bbbc\Gamma)_{<e>} \to  C_{al}^n(\Gamma),~ c \mapsto \tau_c$$
and 
$$C_{al}^n(\Gamma) \to \ov{C}_{\lambda}^n(\bbbc\Gamma)_{<e>}, ~ \tau \mapsto c_{\tau}$$ 
are inverse to each other and compatible with the differentials. The isomorphism of complexes 
$$(C_{al}^*(\Gamma),d_{\Gamma}) \cong (\ov{C}_{\lambda}^*(\bbbc\Gamma)_{<e>},b^t)$$
 induces an injection $H^n(\Gamma) \to \ov{HC}^n(\bbbc\Gamma), ~n \ge 1$. 

The case $n=0$ is different but easy, therefore we leave it to the reader.
  
Let $M$ be a closed Riemannian manifold with fundamental group $\Gamma$ and universal covering $\pi:\tilde M \to M$. Let $C_r^*\Gamma$ be the reduced $C^*$-algebra of $\Gamma$.

We recall the definition of the higher index in $K_0(C_r^*\Gamma)$ of an elliptic differential operator on $M$. 

Let $E$ be $\bbbz/2$-graded hermitian vector bundle on $M$ and let $\tilde E:=\pi^*E$. Then $\tilde E$ is endowed with a right $\Gamma$-action. There is an induced left $\Gamma$-action $R_g^*:\C(\tilde M,\tilde E) \to \C(\tilde M,\tilde E),~(R_g^*s)(x)=R_{g^{-1}}(s(xg))$. Let $D:\C(M,E) \to \C(M,E)$ be an odd symmetric elliptic differential operator. It lifts to a $\Gamma$-invariant odd elliptic operator $\tilde D$ on $\C_c(\tilde M,\tilde E)$. 
On the right $\bbbc \Gamma$-module $\C_c(\tilde M,\tilde E)$ (which is defined using the right $\Gamma$-action $R_{g^{-1}}^*$) we have a $C_r^*\Gamma$-valued scalar product $$\langle x,y \rangle:=\sum_{g \in \Gamma} g\int_{\tilde M} \langle x,R_g^*y \rangle_{\tilde E} dvol_{\tilde M} \ .$$ The completion of $\C_c(\tilde M,\tilde E)$ with respect to the corresponding norm is a Hilbert $C_r^*\Gamma$-module denoted by $H$. The higher index of $D$ is defined as the index of the closure of $\tilde D^+:\C_c(\tilde M,\tilde E^+) \to H^-$. In the following we show that $\tilde D$ is unitarily equivalent to an elliptic differential operator on $M$. In particular its closure is indeed regular and Fredholm (see Appendix \S \ref{pseudodiff}).

If $V$ is a vector space with a left $\Gamma$-action, there is a left $\Gamma$-action on $(\tilde M \times  V)$ defined by $(x,v)\mapsto (xg^{-1},gv)$. 
The Mishenko-Fomenko bundle on $M$ is $${\mathcal P}:=\tilde M \times_{\Gamma} C_r^*\Gamma= (\tilde M \times C_r^*\Gamma)/\Gamma \ .$$ It inherits a $C_r^*\Gamma$-valued scalar product from the standard $C_r^*\Gamma$-valued scalar product on $C_r^*\Gamma$.

Let ${\mathcal P}_{alg}=\tilde M \times_{\Gamma} \bbbc\Gamma \subset {\mathcal P}$. 

\begin{lem} There is an isometric isomorphism between the Hilbert $C_r^*\Gamma$-modules $H$ and $L^2(M,E \ten {\mathcal P})$ inducing an isomorphism between $\C_c(\tilde M,\tilde E)$ and $\C(M,E \ten {\mathcal P}_{alg})$.
\end{lem}

\begin{proof}
A left $\Gamma$-action on  $\C(\tilde M, \tilde E \ten \bbbc \Gamma)$ is defined  by $$L_h (sg)= (R_h^*s) hg$$ for $s \in \C(\tilde M,\tilde E)$ and $g \in \Gamma$. Let $\C(\tilde M, \tilde E \ten \bbbc \Gamma)^{\Gamma}$ be the subspace of $\Gamma$-invariant sections. There is an $C_r^*\Gamma$-valued scalar product on $\C(\tilde M, \tilde E \ten \bbbc \Gamma)^{\Gamma}$ given by
$$\langle x,y \rangle = \int_F \langle x, y \rangle_{\tilde E \ten C_r^*\Gamma} dvol_{\tilde M}$$ where $F \subset \tilde M$ is a fundamental domain. We denote the completion of $\C(\tilde M, \tilde E \ten \bbbc \Gamma)^{\Gamma}$ with respect to the induced norm by $H_1$.

Any $s \in \C(M,E \ten {\mathcal P}_{alg})$ lifts uniquely to an element $\sum_{g \in \Gamma} s_g  g \in \C(\tilde M, \tilde E \ten \bbbc \Gamma)^{\Gamma}$ with $s_g \in  \C_c(\tilde M, \tilde E)$. 
The induced map is an isometric isomorphism
$$\C(M,E \ten {\mathcal P}_{alg}) \cong  \C(\tilde M, \tilde E \ten \bbbc \Gamma)^{\Gamma}$$
hence we get an isometric isomorphism $L^2(M,E \ten {\mathcal P}) \cong H_1$.

It can easily be checked that the isomorphism
$$\C_c(\tilde M,\tilde E) \to \C(\tilde M, \tilde E \ten \bbbc \Gamma)^{\Gamma},~ s \mapsto \sum_{g \in \Gamma}(R_g^*s) g$$ is an isometry as well inducing an isometry $H \cong H_1$.
\end{proof}

By the lemma and its proof $\tilde D:\C_c(\tilde M,\tilde E) \to H$ is unitarily equivalent to an elliptic differential operator ${\mathcal D}:\C(M,E \ten {\mathcal P}_{alg})  \to L^2(M,E \ten {\mathcal P})$. 

In order to apply Th. \ref{indtheor} we embed ${\mathcal P}$ into a trivial bundle as follows:

Let $\{U_i\}_{i \in I}$ be a finite open covering of $M$ such that $\pi^{-1}U_i$ is diffeomorphic to $U_i \times \Gamma$. By refining the covering we may assume that $U_i\cap U_j$ is connected for each $i,j \in I$. Let $\{\chi_i^2\}_{i \in I}$ be a subordinate partition of unity.  For each $i\in I$ fix an open set $U'_i \subset \pi^{-1}U_i$, such that $\pi:U_i' \to U_i$ is a diffeomorphism. The projection $p:\tilde M \times C_r^*\Gamma \to {\mathcal P}$ induces isometric isomorphisms $p_i: U'_i \times C_r^*\Gamma \to {\mathcal P}|_{U_i}$. Hence we get an isometry ${\mathcal P}\to M \times C_r^*\Gamma^{|I|}$ by mapping $s_x \in {\mathcal P}_x$ to $(x,\chi_i(x)v_i)_{i \in I}$, where $v_i$ is defined by the equation $p_i^{-1}s_x=(x',v_i)$ for $x \in U_i$. 

Let $g_{ij} \in \Gamma$ be the deck transformation inducing a diffeomorphism $U_i' \cap \pi^{-1}U_j \to U_j' \cap \pi^{-1}U_i$. One verifies easily that $P=(\chi_i\chi_j g_{ij})_{ij} \in \C(M,M_{|I|}(C_r^*\Gamma))$ is the projection onto the image of the embedding ${\mathcal P}\to M \times C_r^*\Gamma^{|I|}$. 

Since ${\mathcal D}$ and $P (\oplus ^{|I|}D) P$ are elliptic pseudodifferential operators on $L^2(M,E \ten {\mathcal P})$  with the same symbol, we get $$\ind \tilde D^+=\ind {\mathcal D}^+= \ind ((P (\oplus^{|I|}D) P)^+) \in K_0(C_r^*\Gamma) \ .$$
 
Let $\Bi \subset C_r^*\Gamma$ be a smooth subalgebra containing $\bbbc\Gamma$. Let $\tau \in C^n_{al}(\Gamma)$ be a cocycle such that $c_{\tau}$ extends to a continuous cyclic cocycle on $\Bi$. Such a cocycle $\tau$ is called extendable. 

Note that $P \in \C(M) \ten M_{|I|}(\bbbc \Gamma)$, hence its Chern character form is in $\Ob^*(M) \ten \Omega_*\bbbc\Gamma/[\Omega_*\bbbc\Gamma,\Omega_*\bbbc\Gamma]$, which we indicate by writing $\ch_{\bbbc \Gamma}^M(P)$. Here $\Omega_*\bbbc\Gamma$ is defined using the algebraic tensor product. 

By Theorem \ref{indtheor} 
$$c_{\tau} \ch_{\Bi} \ind ((P (\oplus^{|I|}D) P)^+) =  (-1)^k \int_M \Td(M)\pi_! \ch^{TM}[\sigma(D^+)] c_{\tau}\ch^M_{\bbbc\Gamma}(P) $$
with $k=\frac{\dim M(\dim M+1)}{2}$.

It remains to identify $c_{\tau}\ch_{\bbbc\Gamma}^M(P)$.
 
Let $\nu:M \to B\Gamma$ be the classifying map of the covering $\tilde M \to M$.

Lift the functions $\chi_i$ to functions $\chi_i':U_i' \to \bbbr$ and define $h= \sum_i \chi_i'^2$. 

In \cite{lo} it was shown that for $\tau \in C^n_{al}(\Gamma)$ with $d_{\Gamma}\tau=0$ the differential form $$\tilde \omega_{\tau} =\sum_{g_1, \dots, g_n} (R_{g_1}^*d_Mh)\wedge (R^*_{g_2}d_Mh) \wedge \dots (R^*_{g_n}d_Mh) \tau(e,g_1,\dots, g_n)$$ on $\tilde M$ is closed and $\Gamma$-invariant and that the form $\omega_{\tau} \in \Omega^n(M)$ defined by $\pi^*\omega_{\tau}=\tilde \omega_{\tau}$ fulfills
$$\nu^*[\tau]=[\omega_{\tau}] \in H^n(M)$$
under the identification $H^*(\Gamma)\cong H^*(B\Gamma)$.

\begin{prop}
\label{Chern}
Let $\tau \in C^n_{al}(\Gamma) ~n \ge 1$ with $d_{\Gamma}\tau=0$. Then on the level of differential forms $$c_{\tau}(\ch_{\bbbc\Gamma}^M(P))= (- 1)^{\frac{(n-1)n}{2}}((2 \pi i)^n n!)^{-1}~ \omega_{\tau} \ .$$
Hence in $H^n(M)$
$$[c_{\tau}(\ch_{\bbbc\Gamma}^M(P))]=(- 1)^{\frac{(n-1)n}{2}}((2 \pi i)^{n} n!)^{-1}~ \nu^*[\tau] \ .$$
\end{prop}

\begin{proof} The following calculation is a modification of an argument in \cite{lo}.
First we show that the connection $P\di P$ behaves like a flat connection in the sense that  in the expansion of $c_{\tau}(\ch_{\bbbc \Gamma}^M(P))$ with respect to the decomposition $d_{tot} P= d_M P + \di P$ all terms containing the factor $P\di P\di  P$ vanish. 

By the cyclicity of the trace we only need to consider terms of the form $c_{\tau}(\tr P (d_{tot} P)^{m} P(\di P) (\di P))$. Using Leibniz rule for $\di$ we deduce that $\tr P (d_{tot} P)^{m} P(\di P) (\di P)$ can be written as $$\sum_{i_0,i_1, \dots, i_n \in I}\chi_{i_0}\chi_{i_1}f_{i_1 \dots i_{n-1}} \chi_{i_{n-1}}\chi_{i_n}^2 \chi_{i_0} g_{i_0i_1}\di g_{i_1i_2} \dots \di g_{i_{n-1}i_n} \di g_{i_ni_0}  $$ in $\Omega^*M \ten \Oi\bbbc\Gamma/[\Oi\bbbc\Gamma,\Oi\bbbc\Gamma]$ with $f_{i_1 \dots i_{n-1}} \in \Omega^*M$. Since
$$\chi_{i_0}\chi_{i_1}f_{i_1 \dots i_{n-1}} \chi_{i_{n-1}}\chi_{i_n}^2 \chi_{i_0} c_{\tau}(g_{i_0i_1}\di g_{i_1i_2}, \dots, \di g_{i_{n-1}i_n} \di g_{i_ni_0})$$
$$=\chi_{i_1} f_{i_1 \dots i_{n-1}} \chi_{i_{n-1}}\chi_{i_n}^2 \chi_{i_0}^2\tau(e,g_{i_1i_2},g_{i_1i_2} g_{i_2i_3}, \dots, g_{i_1i_2} \dots  g_{i_{n-1}i_n}g_{i_ni_0})$$
$$=\chi_{i_1} f_{i_1 \dots i_{n-1}} \chi_{i_{n-1}}\chi_{i_n}^2 \chi_{i_0}^2 \tau(e,g_{i_1i_2},g_{i_1i_3}, \dots, g_{i_1i_n},g_{i_1i_0})$$
and since $\tau$ is antisymmetric, the summands are antisymmetric with respect to the pair $(i_0,i_n)$, hence $c_{\tau}(\tr P (d_{tot} P)^m (\di P) (\di P))=0$.

Thus  $$c_{\tau}(\ch_{\bbbc\Gamma}^M(P)) =  \frac{(-1)^n}{(2\pi i)^n n!}c_{\tau}\tr \left(P(d_MP)(\di P) + P (\di P)(d_MP)\right)^{n} \ .$$ 

We have that
\begin{eqnarray*}
\lefteqn{(P(d_MP)(\di P)P)_{im}}\\
&=& \sum_{j,k,l \in I}\chi_i \chi_j g_{ij} d_M(\chi_j\chi_k)g_{jk} \chi_k \chi_l (\di g_{kl}) \chi_l \chi_m g_{lm} \\
&=& \sum_{j,k,l \in I} \chi_i \chi_j d_M (\chi_j\chi_k)\chi_k \chi_l^2 \chi_m(g_{ik} \di g_{km}- g_{il} \di g_{lm})\\
&=& \frac 12 \sum_{j,k,l \in I}(\chi_i (d_M \chi_j^2)\chi_k^2 \chi_l^2 \chi_m + \chi_i \chi_j^2  (d_M  \chi_k^2) \chi_l^2 \chi_m)(g_{ik} \di g_{km}- g_{il} \di g_{lm})\\
&=&  \frac 12 \sum_{k \in I}\chi_i  (d_M  \chi_k^2) \chi_m g_{ik} \di g_{km} \ .
\end{eqnarray*}
Here we used $\sum_{i \in I} \chi_i^2=1$.

A similar calculation for $P(\di P)(d_M P)P$ shows that
 $$P(d_MP)(\di P)P=P(\di P)(d_M P)P \ .$$
It follows that 
$$((P(d_MP)(\di P)P+P(\di P)(d_M P)P)^{n})_{im}$$
$$= (- 1)^{\frac{(n-1)n}{2}}\sum_{i_1, \dots, i_n} \chi_i d_M(\chi_{i_1}^2) d_M(\chi_{i_2}^2) \dots d_M(\chi_{i_n}^2)\chi_m g_{ii_1}\di g_{i_1i_2}\di g_{i_2i_3} \dots \di g_{i_nm} \ .$$

Hence $c_{\tau}(\ch_{\bbbc\Gamma}^M(P))$ equals, up to the factor $(- 1)^{\frac{(n+1)n}{2}}((2\pi i)^n n!)^{-1}$,
$$ \sum_{i_0,i_1, \dots, i_n} \chi_{i_0}^2 d_M(\chi_{i_1}^2) d_M(\chi_{i_2}^2) \dots d_M(\chi_{i_n}^2) c_{\tau}(g_{i_0i_1}\di g_{i_1i_2}\di g_{i_2i_3} \dots \di g_{i_ni_0})$$
$$=\sum_{i_0,i_1, \dots, i_n} \chi_{i_0}^2 d_M(\chi_{i_1}^2) d_M(\chi_{i_2}^2) \dots d_M(\chi_{i_n}^2) \tau(g_{i_0i_1},g_{i_0i_2},g_{i_0i_3}, \dots, g_{i_0i_n},e) \ .$$
The $\Gamma$-invariant lift of the sum to $\tilde M$ equals
\begin{eqnarray*}
\lefteqn{\sum_{g_0,g_1, \dots, g_n \in \Gamma} R_{g_0}^* h R_{g_1}^* d_Mh R_{g_2}^* d_M h \dots R_{g_n}^*d_M h \tau(g_1,g_1g_2, \dots, g_1 \dots g_n, e)}\\
&=& (-1)^n \sum_{g_0,g_1, \dots, g_n} R_{g_0}^* (h( R_{g_1}^* d_Mh)(R_{g_1g_2}^* d_M h) \dots (R_{g_1 \dots g_n}^*d_M h)) \tau(e,g_1,g_1g_2, \dots, g_1 \dots g_n)\\
&=&(-1)^n \sum_{g_0,g_1, \dots, g_n} R_{g_0}^* (h(R_{g_1}^* d_Mh)(R_{g_2}^* d_M h) \dots (R_{g_n}^*d_M h)) \tau(e,g_1,g_2, \dots, g_n)\\
&=&(-1)^n  \sum_{g_0} R_{g_0}^* (h \tilde \omega_{\tau})\\
&=& (-1)^n  \tilde \omega_{\tau} \ .
\end{eqnarray*}
\end{proof}

\begin{cor} For any cocycle $\tau \in C^n_{al}(\Gamma),~n \ge 1$ such that $c_{\tau}$ extends to a continuous cocycle on $\ov{C}_n^{\lambda}(\Bi)$ 
$$c_{\tau}\ch_{\Bi} \ind \tilde D^+=(-1)^k((2\pi i)^n n!)^{-1}\int_M \Td(M)\pi_! \ch^{TM}[\sigma(D^+)] \nu^* [\tau] $$
with $k=\frac{\dim M(\dim M+1)}{2} + \frac{n(n-1)}{2}$.
\end{cor}

The formula for the higher index in $K_1(C^*_r\Gamma)$ of an ungraded operator is analogous.

If $D$ is the signature operator on $M$, then the right hand side equals up to normalization the higher signature associated to $\tau$. The Novikov conjecture for extendable $\tau$ follows then from the homotopy invariance of $\ind \tilde D^+$.

In the following we give an example of a smooth subalgebra $\Bi$ with $\bbbc \Gamma \subset \Bi \subset C_r^*\Gamma$. Its construction is typical for the construction of smooth subalgebras from unbounded derivations. Up to minor details the construction is due to Connes-Moscovici \cite{cm} who showed that for $\Gamma$ Gromov-hyperbolic every class in $H^*(\Gamma)$ has a representative that is extendable with respect to this particular $\Bi$, which implies the Novikov conjecture for Gromov-hyperbolic groups.

Let $l$ be a word length function on $\Gamma$. Define an unbounded operator $D_l$ on $l^2(\Gamma)$ with domain $\bbbc\Gamma$ by $D_l g:= l(g)g$.  Let $A$ be the smallest subalgebra of $B(l^2(\Gamma))$ containing $\bbbc\Gamma$ and multiplication by elements of $l^{\infty}(\Gamma)$ . We define an unbounded derivation $$\delta:A \to B(l^2(\Gamma)),~\delta(T):=[D_l,T] \ .$$
Note that in general $\delta (g)$ is not $\Gamma$-invariant, hence $\delta(\bbbc\Gamma)$ is not a subset of $C^*_r(\Gamma)$. 

\begin{lem}
\begin{enumerate}
\item $\delta$ is closable.
\item  $\delta(A)\subset A$.
\end{enumerate}
\end{lem}

\begin{proof}
1) Let $(T_n)_{n \in \bbbn}\subset \dom \delta$ be a sequence converging to zero in $B(l^2(\Gamma))$ and with $\lim_{n \to \infty} [D_l,T_n] = L \in B(l^2(\Gamma))$. For $f \in l^2(\Gamma)$ we have that $[D_l,T_n]f= D_l(T_nf)-T_n D_l(f)$. The second term on the right hand side converges to zero, hence the first converges to $Lf$. Since $T_nf$ converges to $0$ and $D_l$ is closable, it follows that $Lf=0$.

2) It is clear that $\delta$ annulates multiplication operators. From $[D_l,g]h=(l(gh)- l(h))gh$ and $|l(gh)-l(h)| \le l(g)$ it follows that $\delta(g) \in A$.
\end{proof}

Let $\delta_1$ be the closure of $\delta$. For $i \in \bbbn$ define inductively the Banach algebra $\A_i=\dom \delta_i$ with norm $\|a\|_i=\|a\|_{i-1} + \|\delta_i(a)\|_{i-1}$ and the operator $\delta_{i+1}$ as the closure of $\delta$ on $\A_i$.  Let $\Ai$ be the projective limit of the Banach algebras $\A_i$. Then $\Ai$ is a locally $m$-convex Fr\'echet algebra. Let $\B_i$ be the Banach algebra $\A_i \cap C^*_r\Gamma$ and $\Bi=\Ai \cap C^*_r\Gamma$. Then $A \cap C_r^*\Gamma \subset \bbbc\Gamma \subset \Bi$, hence $\Bi$ is dense in $C^*_r\Gamma$ and $\bbbc\Gamma$ is dense in $\Bi$. 

\begin{lem}
The algebras $\B_i,~i \in \bbbn,$ and $\Bi$ are closed under holomorphic functional calculus in $C_r^*\Gamma$. 
\end{lem}

\begin{proof}
We set $\B_0=C_r^*\Gamma$ and show that $\B_{i+1}$ is closed under holomorphic functional calculus in $\B_i$ for each $i \in \bbbn_0$ . 

If $x\in \B_{i+1}$ with $\|1-x\|_i <1$, then $x^{-1}= \sum_{n=0}^{\infty}(1-x)^n \in \B_i$. Since $\delta_{i+1}((1-x)^n) \in \A_i$ and $$\|\delta_{i+1}((1-x)^n)\|_i \le n\|(1-x)\|^{n-1}_i\|\delta_{i+1}(1-x)\|_i \ ,$$ we have that $x^{-1}  \in \dom \delta_{i+1} \cap C_r^*\Gamma =\B_{i+1}$. Since $A$ is dense in $\A_i$, for general $x \in \B_{i+1}$ such that $x^{-1} \in \B_i$ exists there is $y \in A \cap \B_i$ such that $\|1-xy\|_i \le \frac 12$. Hence $x^{-1}= y(xy)^{-1}\in \B_{i+1}$.
\end{proof}

We can also obtain an index theorem for operator $D$ twisted by the Mishenko-Fomenko bundle $\Pj(m)=\tilde M \times_{\Gamma} C^*\Gamma$, where $C^*\Gamma$ is the maximal group algebra of $\Gamma$. This is a refined version of the theorem for $C_r^*\Gamma$ since there is a surjective homomorphism $p:C^*\Gamma \to C_r^*\Gamma$. We set $\B(m)=C^*\Gamma$ and $\B(m)_i=p^{-1}\B_i$. The norm on $\B(m)_i$ is given by $\|a\|_i=\|a\|+ \|p(a)\|_i$, where $\|a\|$ is the norm of $a$ in $C^*\Gamma$. It is straight-forward to check that $\B(m)_i$ is closed under holomorphic functional calculus in $\B(m)$.

In the following we briefly discuss how these results modify in the situation of flat foliated bundles. Our approach is motivated by the approach taken in \cite{ji}.

If not specified, the notation is as before. Let $X$ be a  closed manifold. We assume that $\Gamma$ acts on $X$ from the right by diffeomorphisms.  By pull-back one gets an induced left action on the algebra $\Omega^*(X)$ of continuous differential forms on $X$. 

Thus we get reduced crossed products $\B:=C(X) \times_r \Gamma$ and $\Omega^*(X) \times_r \Gamma$. 

We denote elements of $\Omega^*(X)$ by greek letters in the following and assume them homogeneous where necessary. Let $g,h \in \Gamma$.

Multiplication on $\Omega^*(X) \times_r \Gamma$ is given by the formula
$$(\alpha h)(\beta g)=(\alpha \wedge h^*\beta) hg \ .$$
The algebra $\Omega^*(X) \times_r \Gamma$ acts faithfully on the Hilbert $\Omega^*(X)$-module $\Omega^*(X,l^2(\Gamma))$ by $(\omega g) (\alpha v)=(\omega \wedge g^*\alpha) gv$, where $v \in l^2(\Gamma)$.

A smooth subalgebra of $\Omega^*(X) \times_r \Gamma$ was constructed in \cite{ji} if the group acts isometrically with respect to some Riemannian metric on $X$. The following construction, which is a generalization of the construction above, works in general.

Let $A\subset B(\Omega^*(X,l^2(\Gamma)))$ be the algebra
generated by multiplication operators associated to elements in $l^{\infty}(\Gamma)$ and by the algebraic crossed product $\C(X,\Lambda T^*X) \times_{alg} \Gamma \subset \Omega^*(X) \times_r \Gamma$. One checks that  $\delta(T)=[D_l, T]$ is a derivation on $A$. We denote by $\A$ the closure of $A$ in $B(\Omega^*(X,l^2(\Gamma)))$. The de Rham operator $d_X$ on $X$ defines a derivation on $\C(X,\Lambda T^*X) \times_{alg} \Gamma$ by $d_X(\omega g):=(d_X\omega)g$. Leibniz rule follows from
\begin{align*}
d_X(\alpha h \beta g)&=d_X(\alpha\wedge h^*\beta) hg \\
&=(d_X\alpha \wedge h^*\beta + (-1)^{|\alpha|}\alpha \wedge h^*d_X \beta)hg \\
&=(d_X\alpha)h \beta g + (-1)^{|\alpha|}\alpha h  (d_X \beta)g \ .
\end{align*}

We extend $d_X$ to a derivation on $A$ by letting it commute with multiplication operators coming from elements in $l^{\infty}(\Gamma)$. 

\begin{lem}
The derivations $d_X$ and $\delta$ on $\A$ are closable.
\end{lem}

\begin{proof}
Let $(a_n)_{n \in \bbbn}\subset A$ be a sequence converging to zero in $\A$ and with $\lim_{n \to \infty} d_X(a_n) = a \in \A$. For $f \in \C(X,\Lambda T^*X) \ten \bbbc \Gamma \subset \Omega^*(X,l^2(\Gamma))$
$$af=\lim_{n \to \infty}d_X(a_nf) \ .$$ 

Choose a Riemannian metric on $X$ and let $\Omega_{(2)}^*(X)$ be the $L^2$-completion of $\Omega^*(X)$. The unbounded operator $d_X$ is densely defined and closable on the Hilbert space $\Omega^*_{(2)}(X) \ten l^2(\Gamma)$. Since there is a continuous injection $\Omega^*(X,l^2(\Gamma)) \to \Omega^*_{(2)}(X) \ten l^2(\Gamma)$, the operator $d_X$ is closable on $\Omega^*(X,l^2(\Gamma))$ (considered here as a Banach space) as well. Hence $af=0$.

The proof of the closability of $\delta$ is analogous, using the closability of $D_l$ on $\Omega^*_{(2)}(X) \ten l^2(\Gamma)$.
\end{proof}

Let now $\delta_1, d^X_1$ be the closure of $\delta,d_X$ respectively. We define, generalizing the above construction, for $i \in \bbbn$ inductively the Banach algebra $\A_i=\dom \delta_i\cap \dom d^X_i$ with norm 
$$\|a\|_i=\|a\|_{i-1} + \|\delta_i(a)\|_{i-1}+\|d_i^X(a)\|_{i-1}$$ and the operator $\delta_{i+1}, d^X_{i+1}$ as the closure of $\delta, d_X$, respectively, on $\A_i$.  Then $\A_i$ is closed under holomorphic functional calculus in $\A$.  Again, let $\B_i$ be the Banach algebra $\A_i \cap \B$. We have that $A \cap \B=\C(X)\times_{alg} \Gamma$, which is dense in $\B_i$ for any $i$.

If $\Gamma$ is trivial, then there is a continuous embedding from $\B_i$ to $C^1(X)$. Examples of cyclic cocycles on $C^1(X)$ are the traces $f \mapsto \int_X \alpha \wedge df$ for $\alpha$ a closed $k$ form, $k=\dim X-1$.

We get a $\B$-vector bundle $\Pj^X:=\tilde M \times_{\Gamma} \B$. An isometric embedding $\Pj^X \to M \times \B^{|I|}$ can constructed as above. (Indeed, it holds that $\Pj^X=\Pj \ten_{C_r^*\Gamma}\B$.) Now let $F$ be a hermitian vector bundle on $\tilde M \times_{\Gamma} X$ and $F \to (\tilde M \times_{\Gamma} X)\times \bbbc^p$ an isometric embedding. We let $P_F:(\tilde M\times_{\Gamma} X) \times \bbbc^p \to F$ be the orthonormal projection and denote by $\tilde P_F$ the $\Gamma$-invariant lift to $(\tilde M \times X) \times \bbbc^p$.

\begin{lem} 
Let $T$ be a $\C(X)$-linear $\Gamma$-equivariant differential operator on $\C_c(\tilde M \times X)$. The $T$ descends to a differential operator ${\cal T}$ on $\C(M,\Pj^X)$. The map $T \mapsto {\cal T}$ is compatible with taking sums and products of differential operators.
\end{lem}

\begin{proof}
The operator $T$ induces a $\Gamma$-equivariant operator on $\C_c(\tilde M \times X)\ten \bbbc\Gamma$ and extends to a $\Gamma$-equivariant operator on $\C(\tilde M,C(X)\times_r \Gamma)$, thus is well-defined on $\C(\tilde M,C(X)\times_r \Gamma)^{\Gamma} \cong \C(M,\Pj^X)$. 
\end{proof}

Thus the projection $\tilde P_F$ defines a projection $P_{\F}$ on $C(M,\Pj^X)^p$ and thus a $C(X) \times_r \Gamma$-vector bundle $\F$ on $M$. We denote the composition of the projection $M \times \B^{|I|p}\to (\Pj^X)^p$ with $P_{\F}$ by $P_{\F}$ again. We have that $P_{\F} \in \C(M,M_{|I|p}(\Bi))$. 

Now one can apply Theorem \ref{indtheor} to the pairing of $[D] \in KK_0(C(M),\bbbc)$ with  $[P_{\F}] \in K_0(C(M,C(X) \times_r \Gamma))$, where $D$ is as before.

For cyclic cocycles concentrated at the conjugacy class of the identity (see \cite{golo} for the terminology) on $\Bi$ one can evaluate the pairing with $\ch^M_{\Bi}(P_{\F})$ further, in a similar but more complicated way as in Prop. \ref{Chern}, see \cite{ji}\cite{golo} for related calculations. We refrain from giving details since the formula would be a special case of \cite{golo}. (To be precise, in \cite{golo} it was assumed that the holonomy groupoid is Hausdorff and remarked that the results might hold in general. Here we do not make this assumption.) 

As above we also get index formulas if we take the maximal crossed product.

\subsection{An index theorem for $C^*$-dynamical systems}

Let $G$ be an $n$-dimensional oriented compact Lie group with $n$ odd. We assume that $G$ is endowed with an invariant Riemannian metric with unit volume. Let $D$ be an invariant Dirac operator on an invariant Clifford bundle $E$ on $G$.

Let $(\A,G,\alpha)$ be $C^*$-dynamical system associated to $G$. Hence $\A$ is a $C^*$-algebra, which we assume to be unital, and $\alpha:G \to \Aut \A$ is a group homomorphism such that the map $\A \to C(G,\A),~ a \mapsto \alpha(a)$ is a well-defined homomorphism of $C^*$-algebras. 

The operator $P=1_{\ge 0}(D)$ acts on the Hilbert $\A$-module $H:=L^2(G,E\ten\A)$. We define the following Toeplitz type extension: Let $T(D,\alpha)$ be the $C^*$-algebra generated by the compact operators $K(PH)$ on $PH$ and the Toeplitz operators $P\alpha(a)P,~ a \in \A$. 

Let $\Psi_{\A}(E)$ be the closure in $B(L^2(G,E \ten \A))$ of the algebra of classical pseudodifferential operators of order smaller than or equal to zero. We obtain a commutative diagram with exact rows
$$\begin{CD}
0 @>   >> K(PH) @>  >> T(D,\alpha) @>  >> \A @>  >> 0 \\
  @.    @VV   V@VV   V@VV  V   \\
0 @>  >> K(H)  @>  >> \Psi_{\A}(E) @> \sigma >> C(SG,\End E)\ten \A @>  >> 0 \ .
\end{CD} $$
See Appendix \ref{pseudodiff} for the exactness of the second row.
Here the second vertical map is defined by $P\alpha(a)P \mapsto (P\alpha(a)P + (1-P))$ and the last vertical map is defined by
$a \mapsto \sigma(1-P) +(\alpha(a) \circ p) \sigma(P)$  where $p:SG \to G$ is the projection. Since this map is injective, the last map in the first row, defined as $P\alpha(a)P  \mapsto a$, is well-defined.

The connecting map $K_1(\A) \to K_0(K(PH)) \cong K_0(\A)$ maps $[u] \in K_1(\A)$ with $u \in U_k(\A)$ to $\ind((\oplus^k P)\alpha(u)(\oplus^k P)) \in K_0(\A)$, where $(\oplus^k P)\alpha(u)(\oplus^k P)$ is understood as a Fredholm operator on $(PH)^k$  and $\alpha(u)$ is defined by  applying $\alpha$ componentwise. For notational simplicity we assume that $k=1$ in the following. If $\mu:C(G) \to B(L^2(G,E))$ is the multiplication operator, then $(L^2(G,E),\mu,D)$ is an unbounded Kasparov $(C(G),\bbbc)$-module, which can be paired with $[\alpha(u)] \in K_1(C(G,\A))$. By Prop. \ref{oddind}  
$$\ind(P\alpha(u) P)=[\alpha(u)]\ten_{C(G)} [(L^2(G,E),\mu,D)] \ .$$

In order to apply Theorem \ref{indtheor} define
the algebra $$\Ai=\{x \in \A~|~ (g \mapsto \alpha_g(x)) \in \C(G,\A)\} \ .$$ Endowed with the subspace topology of $\C(G,\A)$ this is a smooth subalgebra of $\A$. Let ${\mathfrak g}=T_eG$ and define $$d_e:\Ai \to {\mathfrak g}^* \ten \Ai  ,~ a \mapsto d_G(\alpha(a))(e) \ .$$ We denote by $[{\mathfrak g}]\in \Lambda^n{\mathfrak g}$ the dual of the volume form of $G$ at $e$. Let $\tau$ be an invariant trace on $\Ai$. Then for $u \in U(\Ai)$ we have that $\alpha(u) \in U(\C(G,\Ai))$. We obtain the following index formula:
\begin{eqnarray*}
\lefteqn{\tau(\ind(P\alpha(u)P))}\\
&=& \tau(\ch_{\Ai}[\alpha (u)][D]) \\
&=&(-1)^{\frac{n(n+1)}{2}} \int_G \hat A(G) \ch(E/S)\tau( \ch_{\Ai}^G \alpha(u))\\
&=& (-1)^{\frac{n(n+1)}{2}}\left(\frac{-1}{2\pi i}\right)^n\frac{(n-1)!}{(2n-1)!}  \langle \tau(u^*d_e u((d_e u^*)(d_e u))^{k-1}), [{\mathfrak g}] \rangle \ .
\end{eqnarray*}

The second equality follows from Theorem \ref{indtheor}.

As an example consider $G=S^1$, $E=S^1 \times \bbbc$ and $D= \frac{1}{i}\frac{d}{dx}$. Let $\delta(u):= \frac{d}{dt}\alpha_t(u)|_{t=0}=d_e(u)[T_eS^1]$. Then we get 
\begin{eqnarray}
\label{dynsys}
\tau(\ind(P\alpha(u)P))=-\frac{1}{2\pi i}\tau(u^*\delta(u)) =\frac{1}{2\pi i}\tau(u\delta(u^*)) \ .
\end{eqnarray}

If furthermore $\A=C(S^1)$ with the $S^1$-action $\alpha$ given by translation and $\tau(1)=1$, then for $u=e^{-2\pi it}$ the formula gives $\tau(\ind(P\alpha(u)P))=1$. Hence in this case the above connecting map $K_1(C(S^1)) \to K_0(C(S^1))$ is given by the Bott periodicity isomorphism.

We relate formula \ref{dynsys} to an index theorem proven in \cite{le}: Assume that the trace is faithful, normal and  $\tau(1)=1$. Let $\A_{\tau}$ be the Hilbert space completion of $\A$ with respect to the scalar product $\langle a,b\rangle:=\tau(a^*b)$. Let $\tilde D$ be the closure of $\frac{1}{i}\frac{d}{dx}$ acting on $L^2(\bbbr,\A_{\tau})$ and let $\tilde P=1_{\ge 0}(\tilde D)$. Let $\tilde \alpha$ be the lift of $\alpha$ to an action of $\bbbr$ on $\A$, furthermore $\pi$ the left regular representation of $\A$ on $\A_{\tau}$ and $\lambda$ the representation of $\bbbr$ on $L^2(\bbbr)$ given by $(\lambda(y)f)(x)=f(x-y)$. Then $(\pi,\lambda)$ is a covariant representation of $(\A,\bbbr,\tilde\alpha)$ and induces a representation $\pi \times \lambda$ of the cross product $\A \times_{\tilde \alpha} \bbbr$ on $L^2(\bbbr,\A_{\tau})$. Let ${\mathcal N}$ be the von Neumann algebra generated by the image of $\pi \times \lambda$ in $B(L^2(\bbbr,\A_{\tau}))$. By \cite{le} the operator $\tilde P \tilde\alpha(u)\tilde P$ is Breuer-Fredholm in $\tilde P{\mathcal N}\tilde P$ and  $$\Ind_{\tau} (\tilde P \tilde\alpha(u)\tilde P)=\frac{1}{2\pi i}\tau(u\delta(u^*)) \ ,$$ where $\ind_{\tau}$ denotes the index with respect to the trace $\tau$.
We conclude:

\begin{prop}
$$\Ind_{\tau} (\tilde P \tilde\alpha(u)\tilde P) = \tau(\ind(P\alpha(u)P)) \ .$$
\end{prop}

\section{Appendix}

\subsection{Index theory and $KK$-theory}
\label{indKK}

Let $\A,\B$ be unital $C^*$-algebras. We recall the notion of a truly unbounded Kasparov $(\A,\B)$-module from \cite{wa}. Truly unbounded Kasparov $(\A,\B)$-modules define elements in $KK_*(\A,\B)$.
We then express the pairing of $K_*(\A)$ with a truly unbounded Kasparov $(\A,\B)$-module in terms of index theory.
We refer to \cite{bl} for more about $KK$-theory.

Let $H$ be a countably generated Hilbert $\B$-module. 

Recall that a densely defined selfadjoint operator $D$ on $H$ is regular if $(1+D^2)$ has a bounded inverse.
For an unbounded selfadjoint regular operator $D$ on $H$ we denote by $H(D)$ the Hilbert $\B$-module whose underlying $\B$-module is $\dom D$ and whose $\B$-valued scalar product is given by $$\langle x,y\rangle_D:=\langle x,y\rangle + \langle Dx,Dy\rangle\ ,$$ where $\langle~,~\rangle$ is the $\B$-valued scalar product on $H$. We say that a regular selfadjoint operator $D$ is Fredholm if $F:=D(1+D^2)^{-\frac 12}$ is invertible in the Calkin algebra $B(H)/K(H)$. This is equivalent to $D:H(D) \to H$ being Fredholm and also to the existence of an odd monotonous smooth function $\chi$ whose limit at $\pm \infty$ is $\pm 1$ such that $\chi(D)^2-1 \in K(H)$. Such a function is called a normalizing function for $D$.

For the definition of the index note that by the Stabilization Theorem  $H \oplus H_{\B} \cong H_{\B}$. Hence there is an inclusion $K(H) \to K(H_{\B})$. The induced map $K_i(K(H)) \to K_i(K(H_{\B}))$ does not depend on the choice of the isomorphism. Thus we get a map $$K_i(B(H)/K(H)) \to K_{i+1}(K(H)) \to K_{i+1}(K(H_{\B})) \cong K_{i+1}(\B) \ .$$
If $H$ is $\bbbz/2$-graded with $H^+ \cong H^-$ and $D$ is odd, then we identify $H^+$ with $H^-$ and define $\ind D^+$  as the image of $[D^+(1+D^2)^{-\frac 12}]\in K_1(B(H^+)/K(H^+))$ in $K_0(\B)$. If $H^+$ is not isomorphic to $H^-$, we define the index of $D$ as the index of the direct sum of $D$ with an invertible odd operator on $H_{\B}^+ \oplus H_{\B}^-$. This works by the Stabilization Theorem.

In the case where $H$ is ungraded, the index $\ind(D) \in K_1(\B)$ is defined as the image of $[2\chi(D)-1]\in K_0(B(H)/K(H))$ in $K_1(\B)$.

\begin{ddd}
Let $H$ be a Hilbert $C^*$-module and $\rho:\A \to B(H)$ a unital $C^*$-homomorphism. 

Let $D$ be a selfadjoint regular Fredholm operator and assume that there is a dense subset $\Ai \subset \A$ such that for all $a \in \Ai$ the operator $[D,\rho(a)]$ is defined on a core for $D$ and extends to a compact operator from $H(D)$ to $H$ and that there is $x \in [0,\frac 12)$ such that $[D,\rho(a)](1+D^2)^{-x}$ is bounded.

Then $(H,\rho,D)$ is called a truly unbounded odd Kasparov $(\A,\B)$-module. 

If in addition $H$ is $\bbbz/2$-graded, $\rho$ is even and $D$ is odd, then $(H,\rho,D)$ is called a truly unbounded even Kasparov $(\A,\B)$-module.

A truly unbounded Kasparov $(\A,\B)$-module $(H,\rho,D)$ is called an unbounded Kasparov $(\A,\B)$-module if $(D^2+1)^{-1} \in K(H)$ and if $[\rho(a),D]$ is bounded for $a \in \Ai$.

A truly unbounded Kasparov $(\A,\B)$-module $(H,\rho,F)$ is a called a bounded Kasparov $(\A,\B)$-module if $F$ is bounded and $F^2-1 \in K(H)$.

If $(H,\rho,D)$ is a truly unbounded odd resp. even Kasparov $(\A,\B)$-module and $\chi$ is a normalizing function for $D$, then $(H,\rho,\chi(D))$ is a bounded odd resp. even Kasparov $(\A,\B)$-module. The class $[(H,\rho,D)]$ in $KK_1(\A,\B)$ resp. $KK_0(\A,\B)$ is defined as the class $[(H,\rho,\chi(D))]$.
\end{ddd}

For example let $D$ be an elliptic scalar selfadjoint pseudodifferential operator on a closed manifold $M$ and let $\A=C(M)$, $\Ai=\C(M)$ and for $a \in \A$ let $\rho(a) \in B(L^2(M))$ be the multiplication operator. Then $(L^2(M),\rho,D)$ is an odd truly unbounded Kasparov $(C(M),\bbbc)$-module.

If $\A=\bbbc$ and $\rho:\bbbc \to B(H)$ is the unique unital homomorphism, then we suppress $\rho$ in the notation in the following. We identify $KK_0(\bbbc,\B)$ with $K_0(\B)$ via the natural isomorphism $[(H,D)] \mapsto \ind(D^+)$ and $KK_1(\bbbc,\B)$ with $K_1(\B)$ via $[(H,D)] \mapsto \ind(D)$.

For the following lemma note that $(H,\rho,D)$ is an even resp. odd truly unbounded Kasparov $(\A,\B)$-module, then $(H^n,M_n(\rho),\oplus^n D)$ is an even resp. odd Kasparov $(M_n(\A),\B)$-module and that a projection $P \in M_n(\A)$ defines a class $[P]_1 \in K_0(\A)$ as well as $[P]_n \in K_0(M_n(\A))$.

In the following we write $a$ for $\rho(a)$. 
 
\begin{lem} Let $P \in M_n(\A)$ be a projection.

Let $(H,\rho,D)$ be an even resp. odd truly unbounded $(\A,\B)$-Kasparov module.

Then in $K_0(\B)$ resp. $K_1(\B)$
$$[P]_n\ten_{M_n(\A)}[(H^n,M_n(\rho),\oplus^n D)]= [P]_1\ten_{\A}[(H,\rho,D)] \ .$$
\end{lem}

\begin{proof} 
By Morita-equivalence the map $i:\A \to M_n(\A),~a \mapsto aE_{11}$ induces an isomorphism in $KK$-theory. We have that 
$$[(H,\rho,D)]=i^*[(H^n,M_n(\rho),\oplus^n D)] \ ,$$ hence
\begin{eqnarray*}
[P]_1\ten_{\A}[(H,\rho,D)] &=& [P]_1\ten_{\A} i^*[(H^n,M_n(\rho),\oplus^n D)] \\
&=& i_*[P]_n \ten_{M_n(\A)}[(H^n,M_n(\rho),\oplus^n D)] \ .
\end{eqnarray*}
Since $i_*[P]_1=[P]_n$
$$[P]_1\ten_{\A}[(H,\rho,D)]=[P]_n\ten_{M_n(\A)}[(H^n,M_n(\rho),\oplus^n D)] \ .$$
\end{proof}

\begin{lem} Let $P \in \A$ be a projection.

Let $(H,\rho,F)$ be an even resp. odd bounded $(\A,\B)$-Kasparov module.

If $F=PFP + (1-P)F(1-P)$, then in $KK_0(\bbbc,\B)$ resp. $KK_1(\bbbc,\B)$
$$[P] \ten_{\A}  [(H,\rho,F)]=[PH,PFP] \ .$$
\end{lem}

\begin{proof}
Define the $C^*$-algebra $\A_P=P\A P \subset \A$ and let $i:\A_P \to \A$ be the injection. Let $p:\bbbc \to \A_P$ the unique unital homomorphism. Then $[(\A_P,p,0)] \in KK_0(\bbbc,\A_P)$ and $i_*[(\A_P,p,0)]=[P]$. Hence 
\begin{eqnarray*}
[P] \ten_{\A}  [(H,\rho,F)] &=& [p] \ten_{\A_P}  i^*[(H,\rho,F)]\\
&=& p^*i^* [PH,\rho,PFP]\\
&=& [PH,\rho \circ i \circ p,PFP]\\
&=& [PH,PFP] \ .
\end{eqnarray*}
\end{proof}

\begin{lem}
Let $A$ be an unbounded symmetric operator on $H$ with $\dom D \subset \dom A$. Assume that $A:H(D) \to H$ is compact. Then $D+A$ is regular. If furthermore there is $x<\frac 12$ such that $A(1+D^2)^{-x}$ is bounded, then $$f(D+A)-f(D) \in K(H_{\A})$$
for any function $f \in C(\bbbr)$ such that $\lim_{x \to \infty}f(x)$ and $\lim_{x \to -\infty}f(x)$ exist. 
\end{lem}

\begin{proof}
Let $(\phi_n)_{n \in \bbbn} \subset C_c(\bbbr)$ be a uniformly bounded sequence converging uniformly to $1$ on each compact subset of $\bbbr$. Then $\phi_n(D)A\phi_n(D)$ converges to $A$ in $B(H(D),H)$. Furthermore $\phi_n(D)A\phi_n(D) \in B(H)$. Since $D$ is regular, the operators $D\pm i$ are invertible. Hence there is $n$ such that $D+A- \phi_n(D)A\phi_n(D)\pm i$ are invertible, thus $D+A- \phi_n(D)A\phi_n(D)$ is regular. A bounded perturbation of a regular operator is regular, thus $D+A$ is regular.

For $D_0=D+A$ and $D_1=D$ let $F_i=D_i(1+D_i^2)^{-1/2}$.
As in the proof of \cite[Prop. 3.7.]{wa} it follows that
$F_0-F_1 \in K(H_{\A}) \ .$
Then $\pi(F_0)=\pi(F_1)$, where $\pi:B(H) \to B(H)/K(H)$ is the projection, hence $\pi(f(F_0))=\pi(f(F_1))$ for any function $f \in C([-1,1])$. 
\end{proof}

\begin{prop} Let $P \in M_n(\Ai)$ be a projection.

Let $(H,\rho,D)$ be an even resp. odd truly unbounded Kasparov $(\A,\B)$-module. Then in $KK_0(\bbbc,\B)$ resp. $KK_1(\bbbc,\B)$
$$[P]\ten_{\A}[(H,\rho,D)]=[PH^n,P(\oplus^n D)P] \ .$$
Hence for $D$ even
$$[P]\ten_{\A}[(H,\rho,D)]=\ind(P(\oplus^n D^+)P) \in K_0(\B) \ ,$$
and for $D$ odd
$$[P]\ten_{\A}[(H,\rho,D)]=\ind(P(\oplus^n D)P) \in K_1(\B) \ .$$
\end{prop}

\begin{proof}
By the first lemma it is enough to consider the case $P\in \A$.

Let $D_P=PDP+(1-P)D(1-P)$.

Since $D_P=D-2P[D,P]$ and $2P[D,P]:H(D) \to H$ is compact, the operator $D_P$ is Fredholm. Furthermore by assumption there is $x< \frac 12$ such that $2P[D,P](1-D^2)^{-x}$ is bounded. Let $\chi$ be a normalizing function of $D$. Then by the previous lemma $\chi(D)-\chi(D_P) \in K(H)$, hence in $KK_0(\A,\B)$ resp. $KK_1(\A,\B)$ we have that $$[(H,\rho,D)]=[(H,\rho,\chi(D))]=[(H,\rho,\chi(D_P)] \ .$$ 
Hence
\begin{eqnarray*}
[P]\ten_{\A}[(H,\rho,\chi(D_P)] &=& [PH,P\chi(D_P)P] \\
&=& [PH,\chi(PD_PP)]\\
&=&[PH,PDP] \ ,
\end{eqnarray*}
where the first equation follows from the second lemma.
\end{proof}

In the following we use the definition of and results on the relative index of projections and the noncommutative spectral flow from \cite{wa}. We also refer to \cite{wa} for history and references concerning the noncommutative spectral flow, which generalizes the family spectral flow introduced by Dai--Zhang \cite{dz}. We denote the even and the odd spectral both by $\spfl$. The relative index of a pair of projections and the relative index of pair of Lagrangian projections are denoted by $\ind$. 

Let $D$ be a regular selfadjoint Fredholm operator on $H$.
Recall that a selfadjoint operator $A \in K(H)$ is called a trivializing operator of $D$  if $D+A$ is invertible. If $H$ is $\bbbz/2$-graded and $D$ is odd, we assume furthermore that $A$ is odd.

\begin{prop} 
\label{oddind}
Let $(H,\rho,D)$ be an unbounded even resp. odd Kasparov $(\A,\B)$-module. 

Let $U \in M_n(\Ai)$ be a unitary such that $U^* \in M_n(\Ai)$ as well and let $[U]$ be its class in $K_1(A)$. 

If there is a trivializing operator $A$ of $\oplus^n D$, 
then with $P=1_{\ge 0}(\oplus^n D+A)$
\begin{eqnarray*}
\ind(P,UPU^*) &=& \spfl((1-t)(\oplus^n D) +tU(\oplus^n D)U^*,A,UAU^*)\\
&=&[U]\ten_{\A}[(H,\rho,D)] \ .
\end{eqnarray*}  

Without the assumption on the existence of trivializing operators formula \ref{spflowodd} below holds in the ungraded case and formula \ref{spfloweven} in the graded case. 
\end{prop}

\begin{proof} The equality $\ind(P,UPU^*)=\spfl((1-t)(\oplus^n D) +tU(\oplus^n D)U^*,A,UAU^*)$ was proven in \cite[Example after Prop. 3.15]{wa} in the ungraded case. The proof of this formula in the graded case is analogous.

For the second equality it is enough to consider the case $n=1$ by the first lemma.

Assume that $H$ is ungraded. 

Let $\chi\in \C(\bbbr)$ be a monotonous function with $\chi(x)=0$ for $x<\frac 13$ and $\chi(x)=1$ for $x >\frac 23$. 

Let $\dira_{S^1}= \frac 1i \frac{d}{dx}$ on $L^2(S^1)$ and let $[\dira_{S^1}] \in KK_1(C(S^1),\bbbc)$ be the corresponding class, where the $C(S^1)$-action on $L^2(S^1)$ is given by multiplication. 

Define the $\A$-vector bundle $L(U):=([0,1] \times \A)/(0,v)\sim (1,Uv)$ on $S^1$ and let $\dira_{L(U)}$ be the Dirac operator $\dira_{S^1}$ twisted by the bundle $L(U)$ with the trivial connection. Then $i\dira_{L(U)} + \chi(t)D+(1-\chi(t))UDU^*$ is densely defined and and its closure is Fredholm on the Hilbert $\B$-module $L^2(S^1,L(U)) \ten_{\rho} H$. We claim that 
\begin{eqnarray}
\label{spflowodd}
[U][D] &=& \ind(i\dira_{L(U)} + \chi(t)D+(1-\chi(t))UDU^*) \in K_0(\A) \ . \qquad
\end{eqnarray}

If $D$ admits a trivializing operator $A \in B(H)$, then 
\begin{eqnarray*}
\spfl((1-t)D +tUDU^*,A,UAU^*)&=& \spfl(\chi(t)D +(1-\chi(t))UDU^*) \\
&=& \ind(i\dira_{L(U)} + \chi(t)D+(1-\chi(t))UDU^*)) \ ,
\end{eqnarray*}
where the last equation follows from  \cite[Prop. 3.15]{wa} and the relative $K$-theoretic index theorem.

Let $\chi_0$ resp. $\chi_2$ be a smooth positive function equal to $1$ on $[0,\frac 13]$ resp. on $[\frac 23,1]$ and equal to $0$ on $[\frac 12,1]$ resp. on $[0,\frac 12]$. Let $\chi_1=\sqrt{1-\sqrt{\chi_0 +\chi_2}}$; hence $\chi_1^2 + (\chi_0+ \chi_2)^2=1$. It is easy to check that the map $M$ from the bundle $L(U)$ to the range of the projection $$P(U)=\left(\begin{array}{cc} \chi_1^2 &  \chi_1 (\chi_0 + \chi_2 U) \\ \chi_1 (\chi_0 + \chi_2 U^*) & (\chi_0 +\chi_2)^2 \end{array}\right)$$ on $S^1 \times \A^{2}$ defined by $$(x,v) \mapsto (x,\chi_1(x) v \oplus (\chi_0(x) v + \chi_2(x) Uv))$$ is an isometric isomorphism. It induces an isometric isomorphism between $P(U)(L^2(S^1,\A^2) \ten_{\rho} H)$ and $L^2(S^1,L(U)) \ten_{\rho} H$, denoted by $M$ as well. On $L^2(S^1,L(U)) \ten_{\rho} H$ the maps $$\dira_{L(U)} -M^{-1} P(U)(\oplus^{2}\dira_{S^1}) P(U) M \in M_{2}(C(S^1,\A))$$ and $\chi(t)D+(1-\chi(t))UDU^* - M^{-1}P(U)(\oplus^2 D) P(U)M$ are bounded. Hence in $K_0(\B)$
$$\ind(-i\dira_{L(U)} + \chi(t)D+(1-\chi(t))UDU^*)= \ind(P(U)(-i \oplus^{2}\dira_{S^1} + \oplus^2 D) P(U)) \ .$$

We have that $[U]\ten_{\A}[(H,\rho,D)]=[P(U)]\ten_{C(S^1,\A)} ([\dira_{S^1}]\ten [(H,\rho,D)])$.

Since the Kasparov product
$[\dira_{S^1}]\ten [(H,\rho,D)] \in KK_0(C(S^1,\A),\B)$ is represented by the odd selfadjoint operator
$$\left(\begin{array}{cc} 0 & -i\dira_{S^1} + D \\ i\dira_{S^1} + D & 0 \end{array}\right)$$ on $L^2(S^1,H^+ \oplus H^-)$, the previous proposition implies that 
$$([P(U)]\ten_{C(S^1,\A)} [\dira_{S^1}])\ten [(H,\rho,D)]=\ind(P(U)((-i \oplus^{2}\dira_{S^1}) + (\oplus^2 D)) P(U)) \ .$$

If $H$ is graded, let $\sigma$ be the grading operator. One proves analogously that
\begin{eqnarray}
\label{spfloweven}
[U][D] &=& \ind(-\sigma\dira_{L(U)}+ i\sigma (\chi(t)D+(1-\chi(t))UDU^*)) \in K_1(\B) \ .\qquad  
\end{eqnarray}

This,  \cite[Remark after Prop. 8.4]{wa} and the relative index theorem imply the assertion.
\end{proof}

By the following argument, which was pointed out to the author by Ryszard Nest, the pairing of $U$ with $(H,\rho,D)$ as in the proposition can always be expressed in terms of a spectral flow: Let $d$ be an invertible unbounded operator on $H_{\B}$ with compact resolvents and such that the ranges of $1_{\ge 0}(d)$ and $1_{\le 0}(d)$ contain a copy of $H_{\B}$. Then the pairing of $[U]$ with $[(H,\rho,D)]$ coincides with the pairing of $[U]$ with the odd Kasparov $(\A,\B)$-module $(H \oplus H \oplus H_{\A}, \rho \oplus 0 \oplus 0, D  \oplus (-D) \oplus d)$. The index of $D \oplus (-D) \oplus d$ in $K_1(\B)$ vanishes, and the operator admits spectral sections. Now apply the proposition.

\subsection{Pseudodifferential operators over $C^*$-algebras}
\label{pseudodiff}

We refer to \cite{mf} for the definition and general facts about
pseudodifferential operators over $C^*$-algebras.

Let $\A$ be a $C^*$-algebra with unit. Let $M$ be a closed Riemannian manifold and let $E$ be an $\A$-vector bundle
over $M$ endowed with an $\A$-valued metric. Endow $\A^n$ with the standard $\A$-valued scalar product and
let $E \to M \times \A^n$ be a smooth isometry (such an isometry always exists). Let $F$ be the complement of $E$ in $M \times \A^n$.

Let $\Delta$ be the scalar Laplacian on $M$. For any $p \in \bbbr$ we define the Sobolev space $H^p(M,\A^n)$ as the completion of $\C(M,\A^n)$ with respect to the norm induced by the $\A$-valued scalar product
$$<f,g>_{H^p}:=<(1+\Delta)^{p/2}f,(1+\Delta)^{p/2}g>_{L^2} \ .$$

We first assume that $E=\A^n$.

Let $P:\C(M,E) \to \C(M,E)$ be a symmetric pseudodifferential operator of order $s \in \bbbr$. The nonsymmetric case can be reduced to the symmetric by considering $\left(\begin{array}{cc} 0 & P^* \\ P & 0 \end{array}\right)$. 
 
In \cite{mf} it was shown that $P: H^{p+s}(M,E) \to H^{p}(M,E)$ is continuous.

\begin{lem}
The operator $P: H^{p+s}(M,E) \to H^p(M,E)$ is adjointable.
\end{lem}

\begin{proof}
It is straightforward to check that $P^T=(1+ \Delta)^{-2s-p}P(1+ \Delta)^p:H^p(M,E) \to H^{p+s}(M,E)$ is the adjoint of $P$.
\end{proof}

The proof of the following lemma is analogous to the classical case and is given here for completeness.

\begin{lem}
Let $s>0$. If $P$ is elliptic, then
$P$ as an unbounded operator on $L^2(M,E)$ with domain $H^s(M,E)$ is selfadjoint.
\end{lem}

\begin{proof}
Let $Q$ be the parametrix of $P$ and $P^*$ the adjoint of $P$. The closure of
$QP^*$ equals the closure of $QP$ on $L^2(M,E)$. Thus $Kh=(1-QP^*)h \in \C(M,E)$ for $h \in \dom P^*$.  Clearly $QP^*h \in \dom P$ . Hence
$$h=QP^*h -Kh \in \dom P \ .$$
\end{proof}

\begin{prop}
Let $s>0$. If $P$ is elliptic, then $P$ is regular as an unbounded operator on $L^2(M,E)$ with domain $H^s(M,E)$.
\end{prop}

\begin{proof}
From $$\|(P^2+1)f\|_{L^2} \ge \|f\|_{L^2},~ f \in \C(M,E)\  ,$$
it follows that the operator $P^2+1:H^{2s}(M,E) \to  L^2(M,E)$ is injective and its range is closed. Furthermore it is adjointable by the first lemma.

 It follows that the range of $(1+P^2):H^{2s}(M,E) \to
L^2(M,E)$ is complemented. By the previous lemma $(1+P^2)$ is selfadjoint, hence
$$\Coker(1+P^2)=\Ker(1+P^2)=\{0\} \ .$$

Therefore $(1+P^2)$ is surjective and thus $P$ is regular. 
\end{proof}

\begin{cor}
Let $s >0$ and assume that $P$ is elliptic. The identity induces an adjointable isomorphism between the Hilbert $\A$-modules $H(P)$ and $H^s(M,E)$. 
\end{cor}

\begin{proof}
The identity $H(P) \to H^s(M,E)$ equals the composition of $(1+P^2)^{\frac 12}:H(P) \to L^2(M,E)$ with $(1+P^2)^{-\frac 12}:L^2(M,E) \to H^s(M,E)$. The first map is an isometry and the second is adjointable since $(1+P^2)^{-\frac 12}$ is a pseudodifferential operator of order $-s$.
\end{proof} 

Since $(1+P^2)^q$ is a pseudodifferential operator for any $q \in \bbbr$ we also get adjointable isomorphisms of Hilbert $\A$-modules $(1+P^2)^q:H^{p+2qs}(M,E) \to H^p(M,E)$. 

We conclude that if $E$ is $\bbbz/2$-graded and $P$ is odd, then the index of $P^+:H^{p+s}(M,E^+) \to H^s(M,E^-)$ is independent of $p$ and equals the index of $P^+(1+P^2)^{-\frac 12}:H^p(M,E^+) \to H^p(M,E^-)$.

For general $E$ let $e \in \C(M,M_n(\A))$ be the orthogonal projection onto $E$. Define $\Delta_e:=e\Delta e +(1-e)\Delta(1-e)$. Since $\Delta_e:H^2(M,\A^n) \to L^2(M,\A^n)$ is regular, the restriction of $\Delta_e$ to $L^2(M,E)$ is regular as well. We define $H^p(M,E)$ as the completion of $\C(M,E)$ with respect to the norm induced by the $\A$-valued scalar product
$$<f,g>_{H^p}:=<(1+\Delta_e)^{p/2}f,(1+\Delta_e)^{p/2}g>_{L^2} \ .$$
After replacing $\Delta$ by $\Delta_e$, the statements of this sections hold for general $E$.

Furthermore from the previous corollary one can deduce that the injection $H^p(M,E) \to H^p(M,\A^n)$ is adjointable.

In the end we note the exactness of the sequence associated to the symbol map for classical pseudodifferential operators over $C^*$-algebras.

We assume that $E$ is the trivial vector bundle with fiber $\A$. The general case can be derived from this. Let $\Psi_{\A}$ be the closure of the algebra of classical pseudodifferential operators of order smaller than or equal to zero in $B(L^2(M,\A))$. Let $SM$ be the sphere bundle of $TM$. Since $C(SM)$ and $K(L^2(M))$ are nuclear, the algebra $\Psi_{\bbbc}$ is nuclear as an extension of $C(SM)$ by $K(L^2(M))$ \cite[Theorem 15.8.2]{bl}. Hence $\Psi_{\bbbc} \ten \A \cong \Psi_{\A}$ where $\ten$ is any tensor product of $C^*$-algebras, and there is the commutative diagram 
$$\begin{CD}
0 @>   >> K(L^2(M,\A)) @>  >> \Psi_{\bbbc} \ten \A @>  >> C(SM,\A) @>  >> 0 \\
  @.    @VV   V@VV   V@VV  V   \\
0 @>  >> K(L^2(M,\A))  @>  >> \Psi_{\A} @> \sigma >> C(SM,\A) @>  >> 0 \ .
\end{CD} $$
In particular the second row is exact.

\textsc{Leibniz-Forschungsstelle der G\"ottinger Akademie der Wissenschaften\\Waterloostr. 8 \\30169 Hannover \\ Germany} 

\textsc{wahlcharlotte@googlemail.com}

\end{document}